\documentclass{llncs}
\usepackage{verbatim}
\usepackage{graphicx}
\usepackage{amsfonts}
\usepackage{amscd}
\usepackage{amssymb}
\usepackage{amsmath}

\usepackage{alltt}
\usepackage{rotating}
\usepackage{floatpag}
 \rotfloatpagestyle{empty}
\usepackage{graphicx}
\usepackage{multind}
\usepackage{times}

\usepackage{verbatim}
\usepackage{latexsym}
\usepackage{crop}
\usepackage{txfonts}
\usepackage[hyphens]{url}
\usepackage{setspace}
\usepackage{ellipsis} % 
\usepackage[mathscr]{euscript} % powerset.
\usepackage{pifont} %ding
\usepackage[displaymath]{lineno}
\usepackage{rotating}

\usepackage{tikz} % Needs pgf version >= 2.10.
\usetikzlibrary{chains,shapes,arrows,shadings,trees,matrix,positioning,intersections,decorations,
  decorations.markings,backgrounds,fit,calc,fadings,decorations.pathreplacing}
\def\op#1{{\hbox{#1}}} 
\def\tc{\hbox{:}}
\newcommand{\ring}[1]{\mathbb{#1}}

\def\true{\text{true}}
\def\false{\text{false}}

\begin{document}

\title{Mathematics in the Age of the Turing Machine}
\author{Thomas C. Hales\thanks{{Research supported in part by 
NSF grant 0804189 and the Benter Foundation.}}}
\institute{University of Pittsburgh\\
\email{hales@pitt.edu}}
\maketitle

\section*{}

\begin{figure}[h!]
  \centering
\includegraphics[scale=0.8]{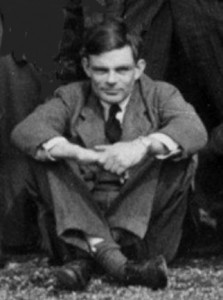}
  \caption{Alan Turing~(image source~\cite{Wel})}
\end{figure}

% image from http://www.heideschwaetzer.de/digitale-fundstuecke/
% beweis-eine-million-dollar-belohnung/
% PERM. obtained from
% Wellcome Library, London.

{

\narrower

\it

``And when it comes to mathematics, 
you must realize that this is the human mind
at the extreme limit of its capacity.'' (H. Robbins) 

\smallskip
\noindent
``\ldots so reduce the use of
the brain and calculate!'' (E. W. Dijkstra)  

\smallskip
\noindent
``The fact that a brain {\it can} do it seems to suggest that the
difficulties [of trying with a machine] 
may not really be so bad as they now
seem.''  (A. Turing)
% Essential Turing page 504, 

}

\newpage

\tikzset{help lines/.style={gray}}
\def\smalldot#1{\draw[fill=black] (#1) node [inner sep=1.3pt,shape=circle,fill=black] {}}
\def\graydot#1{\draw[fill=gray] (#1) node [inner sep=1.3pt,shape=circle,fill=gray] {}}
\def\whitedot#1{\draw[fill=gray] (#1) node [inner sep=1.3pt,shape=circle,fill=white,draw=black] {}}
\tikzset{dartstyle/.style={fill=black,rotate=-90,inner sep=0.7pt,dart,shape border uses incircle}}
\tikzset{grayfatpath/.style={line width=1ex,line cap=round,line join=round,draw=gray}}

\pgfmathsetmacro\cmofpt{(2.54/72.0)}

\def\figFANO{
\begin{tikzpicture}
[scale=2.0]
\coordinate (U) at (0,0);
\coordinate (P1) at (0:1);
\coordinate (P2) at (0:2);
\coordinate (P3) at (60:1);
\coordinate (P4) at (60:2);
\coordinate (P6) at (30:1.73205);
\coordinate (P5) at (30:1.1547);
\draw[thick] (U)--(P2)--(P4)--cycle;
\draw[thick] (U)--(P6);
\draw[thick] (P1)--(P4);
\draw[thick] (P2)--(P3);
\draw[thick] (P5) circle (0.577);
\foreach \i in {U,P1,P2,P3,P4,P5,P6} {\smalldot {\i}; }
\end{tikzpicture}
}

\def\figKISSING{
\begin{tikzpicture}
[scale=1.0]
\coordinate (P4) at (0,0);
\coordinate (P1) at (-150:1);
\coordinate (P2) at (-90:1);
\coordinate (P3) at (-30:1);
\coordinate (P5) at (30:1);
\coordinate (P6) at (90:1);
\coordinate (P7) at (150:1);
\foreach \i in {P1,P2,P3,P4,P5,P6,P7} { \draw[thick,fill=gray!60] (\i) circle (0.5); }
\end{tikzpicture}
}

\def\figEE{
\begin{tikzpicture}
[scale=1.4]
\coordinate (P1) at (0:1);
\coordinate (P2) at (0:2);
\coordinate (P3) at (0:3);
\coordinate (P4) at (0:4);
\coordinate (P5) at (0:5);
\coordinate (P6) at (0:6);
\coordinate (P7) at (0:7);
\coordinate (P8) at (3,1);
\draw[thick] (P1)--(P7);
\draw[thick] (P3)--(P8);
\foreach \i in {P1,P2,P3,P4,P5,P6,P7,P8} { \smalldot {\i}; }
\end{tikzpicture}
}

\def\figPSEUDO{
\begin{tikzpicture}
[scale=0.8]
\coordinate (P1) at (1.6,3.31);
\coordinate (P2) at (2.7,3.4);
\coordinate (P3) at (3.12,3.21);
\coordinate (P4) at (2.0,3.13);
\coordinate (P5) at (0.8,2.33);
\coordinate (P6) at (1.25,2.2);
\coordinate (P7) at (0.8,1.16);
\coordinate (P8) at (1.3,0.97);
\coordinate (P9) at (1.65,0.25);
\coordinate (P10) at (2.75,0.2);
\coordinate (P11) at (3.22,0.4);
\coordinate (P12) at (4.0,1.2);
\coordinate (P13) at (3.95,2.47);
\coordinate (P14) at (3.5,2.26);
\coordinate (P15) at (2.33,2.13);
\coordinate (P16) at (2.4,0.9);
\coordinate (P17) at (3.5,1.02);
\draw[thick,line join=bevel,fill=gray!0.2] (P1)--(P4)--(P3)--(P2)--cycle;
\draw[thick,line join=bevel,fill=gray!0.1] (P1)--(P5)--(P6)--(P4)--cycle;
\draw[thick,line join=bevel,fill=gray!0.2] (P4)--(P6)--(P15)--cycle;
\draw[thick,line join=bevel,fill=gray!0.3] (P4)--(P15)--(P14)--(P3)--cycle;
\draw[thick,line join=bevel,fill=gray!0.5] (P3)--(P14)--(P13)--cycle;
\draw[thick,line join=bevel,fill=gray!0.2] (P5)--(P7)--(P8)--(P6)--cycle;
\draw[thick,line join=bevel,fill=gray!0.3] (P6)--(P8)--(P16)--(P15)--cycle;
\draw[thick,line join=bevel,fill=gray!0.5] (P15)--(P16)--(P17)--(P14)--cycle;
\draw[thick,line join=bevel,fill=gray!0.6] (P14)--(P17)--(P12)--(P13)--cycle;
\draw[thick,line join=bevel,fill=gray!0.4] (P7)--(P9)--(P8)--cycle;
\draw[thick,line join=bevel,fill=gray!0.5] (P8)--(P9)--(P10)--(P16)--cycle;
\draw[thick,line join=bevel,fill=gray!0.6] (P16)--(P10)--(P17)--cycle;
\draw[thick,line join=bevel,fill=gray!0.7] (P17)--(P10)--(P11)--(P12)--cycle;
\end{tikzpicture}
}

\section{Computer Calculation}

\subsection{a panorama of the status quo}

Where stands the mathematical endeavor?

In 2012, many mathematical utilities are reaching consolidation.  It
is an age of large aggregates and large repositories of mathematics:
the arXiv, Math Reviews, and euDML, which promises to aggregate the
many European archives such as Zentralblatt Math and Numdam.  Sage
aggregates dozens of mathematically oriented computer
programs under a single Python-scripted front-end.

% check dozens, hundreds??
% cite euDML http://www.eudml.eu

Book sales in the U.S. have been dropping  for the past several years. 
% \cite{BK11} ``Sales of print books in the U.S. peaked in 2005 and
% have been in steady decline since.''
Instead, online sources such as  Wikipedia and Math Overflow
%Planet
%  Math}, Weinstein's {\it Math World}, 
%and Sloane's {\it Online
%  Encyclopedia of integer sequences} -- 
are rapidly becoming students'
preferred math references. The Polymath blog organizes massive
mathematical collaborations.  Other blogs organize previously isolated
researchers into new fields of research.  The slow, methodical
deliberations of referees in the old school are giving way;
% like a body with high inertia, no external energy source, and heavy
% friction.  Meanwhile,
now in a single stroke,  Tao blogs, gets feedback, and publishes.

Machine Learning is in its ascendancy.  {\it LogAnswer} and {\it
  Wolfram Alpha} answer our elementary questions about the
quantitative world; {\it Watson} our {\it Jeopardy} questions. {\it
  Google Page} ranks our searches by calculating the largest
eigenvalue of the largest matrix the world has ever known.  {\it Deep
  Blue} plays our chess games. The million-dollar-prize-winning {\it
  Pragmatic Chaos} algorithm enhances our {\it Netflix searches}.  The
major proof assistants now contain tens of thousands of formal proofs
that are being mined for hints about how to prove the next generation
of theorems.

Mathematical models and algorithms rule the quantitative world.
%``Without [applied mathematics] we would have no internet, no scanner
%at store counters, no airplane flights or space exploration, no global
%warming measurements, and of course no baseball
%statistics''~\cite{Pi2011}.  
Without applied mathematics, we would be bereft of Shor's
factorization algorithm for quantum computers, Yang-Mills theories of
strong interactions in physics, invisibility cloaks, Radon transforms
for medical imaging, models of epidemiology, risk analysis in
insurance, stochastic pricing models of financial derivatives, RSA
encryption of sensitive data, Navier-Stokes modeling of fluids, and
models of climate change.
%As Pitici reminds us, without applied math, ``no
%internet, no scanner at store counters, \dots no space exploration
%\dots and of course no baseball statistics''~\cite{Pi11}.  
Without it, entire fields of engineering from Control Engineering to
Operations Research would close their doors.  The early icon of
mathematical computing, Von Neumann, divided his final years between
meteorology and hydrogen bomb calculations.
% Reference Birds and Frogs essay of Dyson.
Today, applications fuel the
economy: in 2011 rankings, the first five of the ``10 best jobs'' are
math or computer related: software engineer, mathematician, actuary,
statistician, and computer systems analyst~\cite{CC11}.  
%\footnote{Yet
%  the popular media still portrays the mathematician as one who writes
%  formulas on windows and mirrors in a crazed frenzy, rather than as a
%  driver of the global economy.}

% CC11 http://www.careercast.com/jobs-rated/10-best-jobs-2011

% http://money.cnn.com/magazines/moneymag/bestjobs/2010/snapshots/1.html

% LogAnswer
% http://www.uni-koblenz.de/~bpelzer/publications/FGHP08_IJCAR08_prel.pdf
% http://en.wikipedia.org/wiki/Netflix_Prize

\bigskip

Computers have rapidly become so pervasive in mathematics that future
generations may look back to this day as a golden dawn.  A
comprehensive survey is out of the question.  It would almost be like
asking for a summary of applications of symmetry to mathematics.
Computability -- like symmetry -- is a wonderful structural property
that some mathematical objects possess that makes answers flow more
readily wherever it is found.  This section gives many examples that
give a composite picture of computers in mathematical research,
showing that computers are neither the panacea that the public at
large might imagine, nor the evil that the mathematical purist might
fear.  I have deliberately selected many examples from pure
mathematics, partly because of my own background and partly to correct
the conventional wisdom that couples computers with applied
mathematics and blackboards with pure mathematics.

\subsection{Birch and Swinnerton-Dyer conjecture}

I believe that the Birch and Swinnerton-Dyer conjecture is the deepest
conjecture ever to be formulated with the help of a computer~\cite{BSD}.  The
Clay Institute has offered a one-million dollar prize to anyone who
settles it.

Let $E$ be an elliptic curve defined by an equation $y^2 = x^3 + a x +
b$ over the field of rational numbers.  Motivated by related
quantities in Siegel's work on quadratic forms, Birch and
Swinnerton-Dyer set out to estimate the quantity
\begin{equation}\label{eqn:np}
\prod {N_p/p},
\end{equation}
where $N_p$ is the number of rational points on $E$ modulo $p$, and
the product extends over primes $p\le P$~\cite{Bir}.  Performing
experiments on the EDSAC II computer at the Computer laboratory at
Cambridge University during the years 1958--1962, they observed that
as $P$ increases, the products (\ref{eqn:np}) grow asymptotically in
$P$ as
\[
c(E) \log^r P,
\]
for some constant $c$, where $r$ is the Mordell-Weil rank of $E$; that
is, the maximum number of independent points of infinite order in the
group $E(\ring{Q})$ of rational points.  Following the suggestions of
Cassels and Davenport, they reformulated this numerical asymptotic law
in terms of the zeta function $L(E,s)$ of the elliptic curve.  Thanks
to the work of Wiles and subsequent extensions of that work, it is
known that $L(E,s)$ is an entire function of the complex variable $s$.
The Birch and Swinnerton-Dyer conjecture asserts that the rank $r$ of
an elliptic curve over $\ring{Q}$ is equal to the order of the zero of
$L(E,s)$ at $s=1$.

A major (computer-free) recent theorem establishes that the Birch and
Swinnerton-Dyer conjecture holds for a positive proportion of all
elliptic curves over $\ring{Q}$~\cite{BS:2010}.  This result, although
truly spectacular, is mildly misleading in the sense that the elliptic
curves of high rank rarely occur but pose the greatest difficulties.

\subsection{Sato-Tate}

% Math content checked Sep 2, 2011.

\begin{figure}[h!]
  \centering
\includegraphics[scale=0.33]{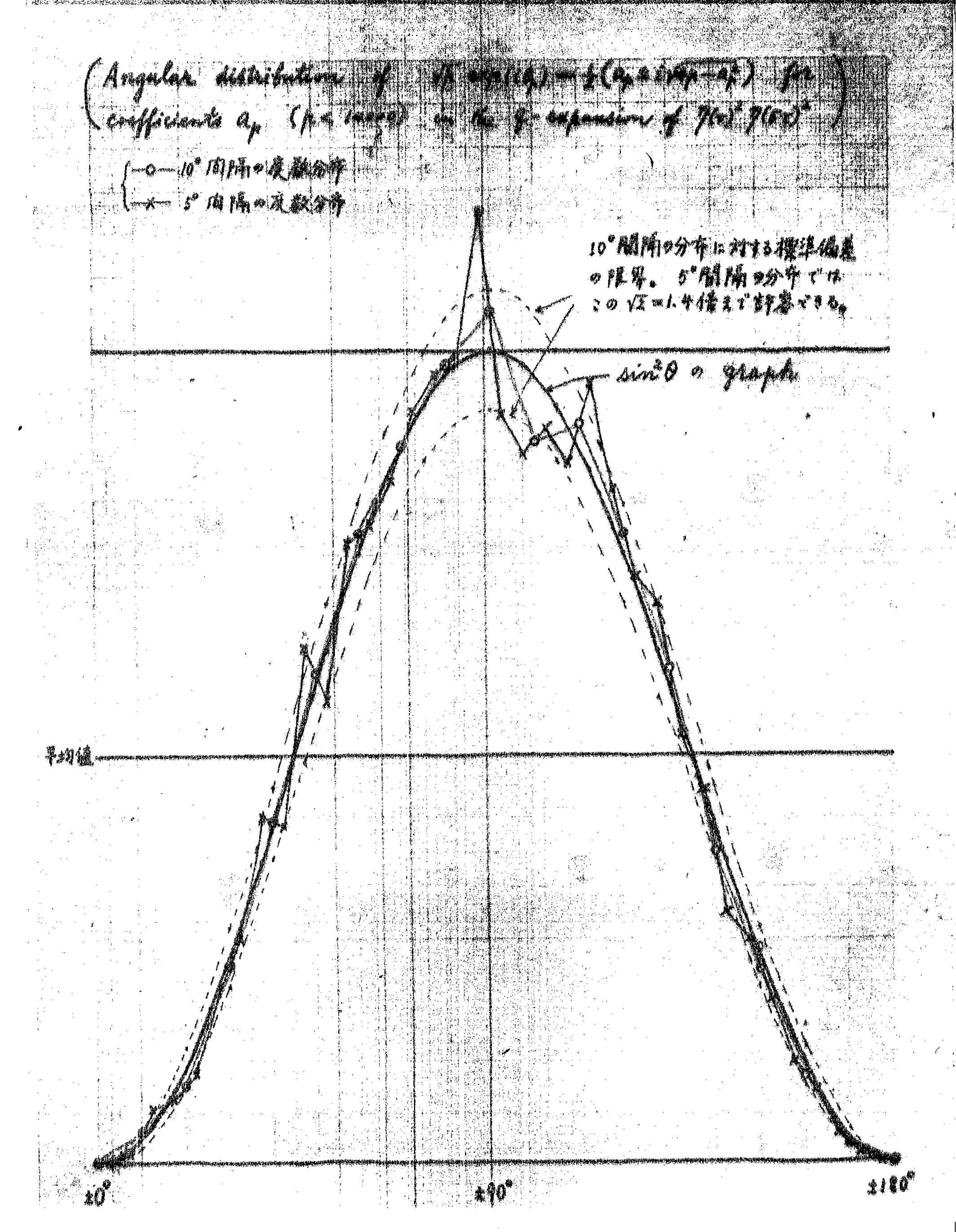}
  \caption{Data leading to the Sato-Tate conjecture (image source~\cite{ST2})}
\label{fig:st}
\end{figure}

%
% Permission and photocopies obtained from Kanji Namba, Jan 2013.
% 463-3 Kitamizote, Sojya Okayama, 719-1117, Japan.
% Sato's calc: graphic, http://www.math.ou.edu/~rschmidt/satotate/ST2.pdf
% Sato's history:  http://www.math.ou.edu/~rschmidt/satotate/page5.html
% accessed 2011

% Hitachi completed their first transistorized electronic computer in 1959.
% http://www.computerhistory.org/brochures/companies.php?alpha=g-i&company=com-42bc1c2c2be73

The Sato-Tate conjecture is another major conjecture about elliptic
curves that was discovered by computer.  If $E$ is an elliptic curve
with rational coefficients
\[
y^2 = x^3 + a x + b,
\]
then the number of solutions modulo a prime number $p$ (including the
point at infinity) has the form
\[
1 + p - 2\sqrt{p}\cos\theta_p.
\]
for some real number $0\le \theta_p\le \pi$.  In 1962, Sato,
Nagashima, and Namba made calculations  of $\theta_p$ on a Hitachi
HIPAC 103 computer to understand how these numbers are distributed as
$p$ varies for a fixed elliptic curve $E$~\cite{Sch}.  By the spring of 1963, the
evidence suggested  $\sin^2\theta$ as a good fit of the data (Figure~\ref{fig:st}).
That is, if $P(n)$ is the set of the first $n$ primes, and
$f:[0,\pi]\to\ring{R}$ is any smooth test function, then for large
$n$,
\[
\frac{1}{n}\sum_{p\in P(n)} f(\theta_p) \quad\hbox{ tends to }\quad
\frac{2}{\pi}\int_0^\pi f(\theta)\,\sin^2\theta\,d\theta.
\]
The Sato-Tate conjecture (1963) predicts that this same distribution is
obtained, no matter the elliptic curve, provided the curve does not
have complex multiplication.  Tate, who arrived at the conjecture
independently, did so without computer calculations.

Serre interpreted Sato-Tate as a generalization of Dirichlet's theorem
on primes in arithmetic progression, and gave a proof strategy of
generalizing the analytic properties of $L$-functions used in the
proof of Dirichlet's theorem~\cite{Se68}.  Indeed, a complete proof of
Sato-Tate conjecture has now been found and is based on extremely
deep analytic properties of $L$-functions~\cite{Car:Bourbaki}.  The
proof of the Sato-Tate conjecture and its generalizations has been one
of the most significant recent advances in number theory.

\subsection{transient uses of computers}

It has become common for problems in mathematics to be first verified by
computer and later confirmed without them.  Some examples are
the construction of sporadic groups, counterexamples to a conjecture
of Euler, the proof of the Catalan conjecture, and the discovery of a formula
for the binary digits of $\pi$.

Perhaps the best known example is the construction of sporadic groups as
part of the monumental classification of finite simple groups.  The
sporadic groups are the $26$ finite simple groups that do not fall into
natural infinite families.  For example,  Lyons (1972) predicted
the existence of a sporadic group of order
\[
2^ 8\cdot 3^7\cdot 5^6\cdot  7\cdot 11 \cdot 31 \cdot 37 \cdot 67.
\]
In 1973, Sims proved the existence of this group in a long unpublished
manuscript that relied on many specialized computer programs.  By 1999
, the calculations had become standardized in group theory packages,
such as GAP and Magma~\cite{HS99}.  Eventually, computer-free
existence and uniqueness proofs were found~\cite{MParker},
\cite{AS92}.

%%
% TR0416 (Sims 1999), for the 1973 calculations
% http://dimacs.rutgers.edu/~havas/TR0416.pdf (Havas and Sims, 1999)
% http://en.wikipedia.org/wiki/Lyons_group

Another problem in finite group theory with a computational slant is
the inverse Galois problem: is every subgroup of the symmetric group
$S_n$ the Galois group of a polynomial of degree $n$ with rational
coefficients?  In the 1980s Malle and Matzat used computers to realize
many groups as Galois groups~\cite{MM}, but with an infinite list of
finite groups to choose from, non-computational ideas have been more
fruitful, such as Hilbert irreducibility, rigidity, and automorphic
representations~\cite{KLS}.

\smallskip

Euler conjectured (1769) that a fourth power cannot be the sum of
three positive fourth powers, that a fifth power cannot be the sum of
four positive fifth powers, and so forth.  In 1966, a computer search
\cite{LP66} on a CDC 6600 mainframe uncovered a counterexample
\[
27^5 + 84^5 + 110^5 + 133^5 = 144^5,
\]
which can be checked by hand (I dare you).  The two-sentence
announcement of this counterexample qualifies as one of the shortest
mathematical publications of all times.  Twenty years later, a more
subtle computer search gave another counterexample~\cite{Elkies88}:
\[
2682440^4 + 15365639^4 + 18796760^4 = 20615673^4.
\]

%% 
% http://en.wikipedia.org/wiki/Euler's_sum_of_powers_conjecture

\smallskip

The Catalan conjecture (1844) asserts that the only solution to the equation
\[
x^m - y^n = 1,
\]
in positive integers $x,y,m,n$ with exponents $m,n$ greater than $1$
is the obvious
\[
3^2 - 2^3 = 1.
\]
That is, $8$ and $9$ are the only consecutive positive perfect powers.
By the late 1970s, Baker's methods in diophantine analysis had reduced
the problem to an astronomically large and hopelessly infeasible finite
computer search.  Mih\u ailescu's proof (2002) of the Catalan
conjecture made light use of computers (a one-minute calculation), and
later the computer calculations were entirely eliminated~\cite{Mih},~\cite{TM03}.

\smallskip

%% 
% computer use described in Bulletin article TM03.
% computer eliminated, see Steiger article.

Bailey, Borwein, and Plouffe found an algorithm for calculating the
$n$th binary digit of $\pi$ directly: it jumps straight to the $n$th
digit without first calculating any of the earlier digits.  They
understood that to design such an algorithm, they would need an
infinite series for $\pi$ in which powers of $2$ controlled the
denominators.  They did not know of any such formula, and made a
computer search (using the PSLQ lattice reduction algorithm) for any
series of the desired form.  Their search unearthed a numerical
identity
\[
\pi = \sum_{n=0}^\infty 
\left(
\frac{4}{8n+1} - \frac{2}{8n+4} - \frac{1}{8n+5} - \frac{1}{8n+6}
\right) 
\left(\frac{1}{16}\right)^n,
\]
which was then rigorously proved and used to implement their
binary-digits algorithm.

%%
% needs citation.
% make the "denominators contained powers of 2" more precise.
% http://www.andrews.edu/~calkins/physics/Miracle.pdf

\subsection{Rogers-Ramanujan identities}

The famous Rogers-Ramanujan identities
\[
1 + \sum_{k=1}^\infty \frac{q^{k^2+a k}}{(1-q)(1-q^2)\cdots (1-q^k)} = 
\prod_{j=0}^\infty \frac{1}{(1-q^{5j+a+1})(1- q^{5j - a +4})},
\qquad a = 0,1.
\]
can now be proved by an almost entirely mechanical procedure from
Jacobi's triple product identity and the $q$-WZ algorithm of Wilf and
Zeilberger that checks identities of $q$-hypergeometric finite
sums~\cite{PP94}.   Knuth's foreword to a book
on the WZ method opens, ``Science is what we understand well enough to
explain to a computer. Art is everything else we do.'' Through the WZ
method, many summation identities have become a science~\cite{PWZ}.

\subsection{packing  tetrahedra}

\begin{figure}[h!]
  \centering
\includegraphics[scale=0.3]{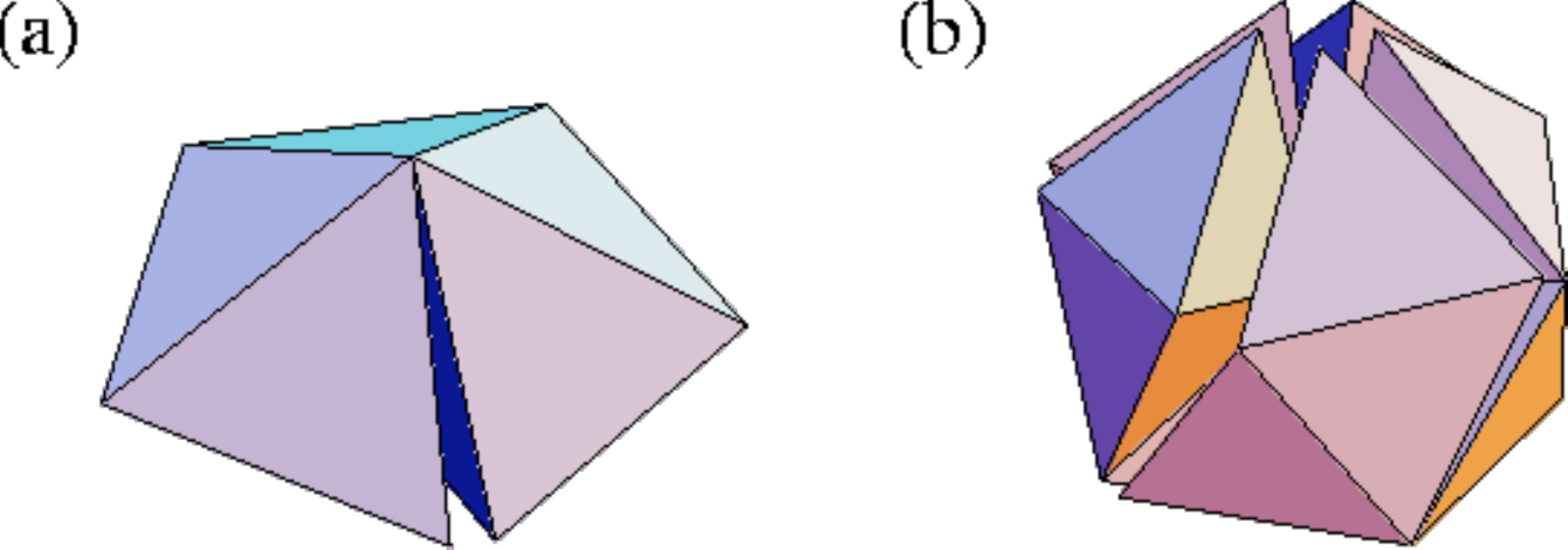}
  \caption{Regular tetrahedra fail to tile space (image source \cite{tetrahedra}).}
\label{fig:gap}
\end{figure}

Aristotle erroneously believed that regular tetrahedra tile space:
``It is agreed that there are only three plane figures which can fill
a space, the triangle, the square, and the hexagon, and only two
solids, the pyramid and the cube''~\cite{Aristotle}.  However,
centuries later, when the dihedral angle of the regular tetrahedron
was calculated:
\[
\arccos(1/3) \approx 1.23 < 1.25664 \approx 2\pi/5,
\]
it was realized that a small gap is left when five regular tetrahedra
are grouped around a common edge (Figure~\ref{fig:gap}).  In 1900, in
his famous list of problems, Hilbert asked ``How can one arrange most
densely in space an infinite number of equal solids of given form,
e.g., spheres with given radii or regular tetrahedra \dots?''

%%
% PERM. obtained by Melissa 12/2012.

\begin{figure}[h!]
  \centering %{CEG_tetrahedra.pdf}
\includegraphics[scale=0.1]{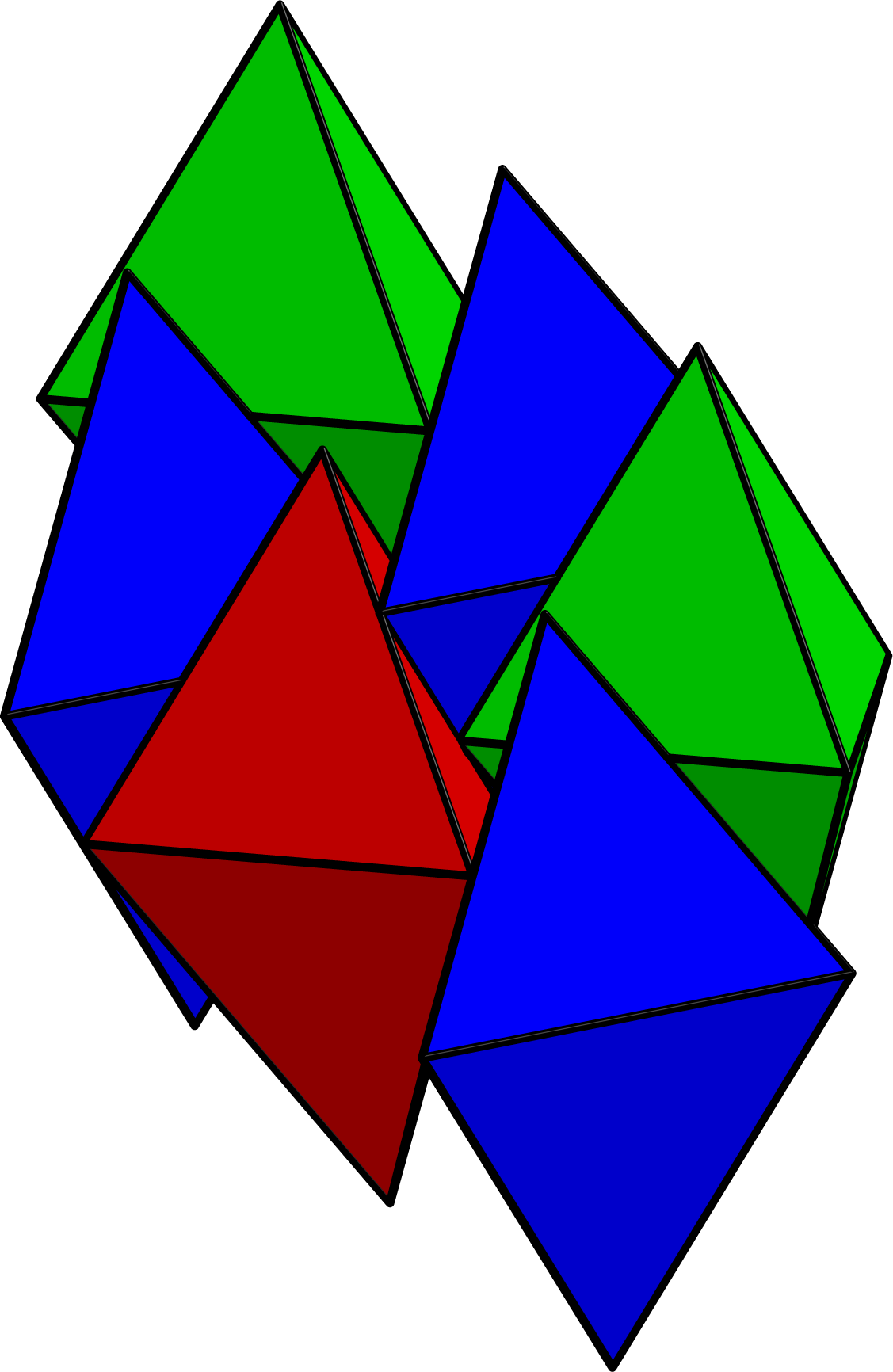}
  \caption{The best packing of
tetrahedra is believed to be the Chen-Engel-Glotzer arrangement with density
$4000/4671\approx 0.856$ (image source~\cite{Chen-2010}).}
\label{fig:CEG}
\end{figure}

Aristotle notwithstanding, until recently, no arrangements of regular
tetrahedra with high density were known to exist.  In 2000, Betke and
Henk developed an efficient computer algorithm to find the densest
lattice packing of a general convex body~\cite{BH2000}.  This opened
the door to experimentation~\cite{Conway-2006}.  For example,
the algorithm can determine the best lattice packing of the convex hull
of the cluster of tetrahedra in Figure~\ref{fig:gap}.  In rapid succession
came new record-breaking arrangements of tetrahedra, culminating in
what is now conjectured to be the best possible~\cite{Chen-2010}.
(See Figure~\ref{fig:CEG}.)  Although Chen had the panache to hand out
Dungeons and Dragons tetrahedral dice to the audience for a hands-on
modeling session during her thesis defense, the best arrangement was
found using Monte Carlo experiments.  In the numerical simulations, a
finite number of tetrahedra are randomly placed in a box of variable
shape.  The tetrahedra are jiggled as the box slowly shrinks until
no further improvement is possible.  Now that a precise conjecture has
been formulated, the hardest part still remains: to give a proof.

% PERM. obtained by Melissa 12/2012.

\subsection{the Kepler conjecture}

\begin{figure}[h!]
  \centering
\includegraphics[scale=0.5]{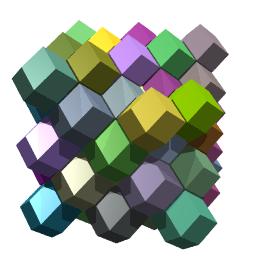}
\caption{An optimal sphere packing is obtained by placing one sphere
  in each three-dimensional honeycomb cell (image
  source~\cite{rhombic}).}
\label{fig:rhombic}
\end{figure}

Hilbert's 18th problem asks to find dense packings of both spheres and
regular tetrahedra.  The problem of determining the best sphere
packing in three dimensions is the Kepler conjecture.  Kepler was led
to the idea of density as an organizing principle in nature by
observing the tightly packed seeds in a pomegranate.  Reflecting on
the hexagonal symmetry of snowflakes and honeycombs, by capping each
honeycomb cell with a lid of the same shape as the base of the cell,
he constructed a closed twelve-sided cell that tiles space.  Kepler
observed that the familiar pyramidal cannonball arrangement is
obtained when a sphere is placed in each capped honeycomb cell
(Figure~\ref{fig:rhombic}).  This he believed to be the densest
packing.

%% wiki source
% PERM. not needed, wiki commons.
% graphic source: http://archinstitute.blogspot.com/2005_10_30_archive.html
% http://upload.wikimedia.org/wikipedia/commons/8/86/Rhombic_dodecahedra.jpg

L. Fejes T\'oth proposed a strategy to prove Kepler's conjecture in
the 1950s, and later he suggested that computers might be used.  The
proof, finally obtained by Ferguson and me in 1998, is one of the
most difficult nonlinear optimization problems ever rigorously solved
by computer~\cite{Hales:2005:Annals}.  The computers calculations
originally took about $2000$ hours to run on Sparc workstations.
Recent simplifications in the proof have reduced the runtime to about
$20$ hours and have reduced the amount of customized code by a factor
of more than $10$.

%%
% Estimates of the amount of custom code in the original 1998 proof
% are over 100 thousand lines.  Sean (in revision paper) say's 137,000
% for Sam and 50,000 for me.  Duplication included.  Now there are
% fewer than $10,000$ lines.

\subsection{the four-color theorem}

The four-color theorem is the most celebrated computer proof in the
history of mathematics.  The problem asserts that it is possible to
color the countries of any map with at most four colors in such a way
that contiguous countries receive different colors.  The proof of this
theorem required about $1200$ hours on an IBM 370-168 in 1976. So much
has been written about Appel and Haken's computer solution to this
problem that it is pointless to repeat it here~\cite{AH4CT}.  Let it
suffice to cite a popular account ~\cite{Wil4CT}, a sociological
perspective~\cite{Mac}, the second generation
proof~\cite{Robertson:1997:JCTB}, and the culminating formal
verification~\cite{gonthier:2008:formal}.

\subsection{projective planes}

A finite projective plane of order $n>1$ is defined to be a set of
$n^2 + n + 1$ lines and $n^2 + n+ 1$ points with the following
properties:
\begin{enumerate}
\item Every line contains $n+1$ points;
\item Every point is on $n+1$ lines;
\item Every two distinct lines have exactly one point of intersection;
\item Every two distinct points lie on exactly one line.
\end{enumerate}

\begin{figure}[h!]
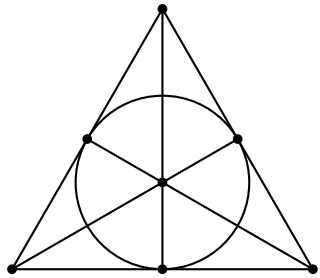

  \centering
\figFANO
\caption{The Fano plane is a finite projective plane of order $2$.}
%  (image source~\cite{fano}).}
\label{fig:P2}
\end{figure}
%%
% Dec 3, 2012, replaced by Figure drawn by Hales.
% Wiki source.
% http://www.log24.com/log/pix11A/110505-WikipediaFanoPlane.jpg = 
% http://en.wikipedia.org/wiki/File:Fano_plane.svg

The definition is an abstraction of properties that evidently hold for
$\ring{P}^2(\ring{F}_q)$, the projective plane over a finite field
$\ring{F}_q$, with $q=n$, for any prime power $q$.  In particular, a
finite projective plane exists whenever $n$ is a positive power of a
prime number (Figure~\ref{fig:P2}).

The conjecture is that every finite projective plane
of order $n>1$ is a prime power.  The smallest integers $n>1$
that are {\it not} prime powers are
\[
6,~10,~12,~14,~15,~\dots
\]
The brute force approach to this conjecture is to eliminate each of
these possibilities in turn.  The case $n=6$ was settled in 1938.
Building on a number of theoretical advances~\cite{MST}, Lam
eliminated the case $n=10$ in 1989, in one of the most difficult
computer proofs in history~\cite{Lam89}.  This calculation was
executed over a period of years on multiple machines and eventually
totaled about 2000 hours of Cray-1A time.

Unlike the computer proof of the four-color theorem, the
projective plane proof has never received independent verification.
Because of the possibilities of programming errors and soft errors
(see Section~\ref{sec:soft}), Lam is unwilling to call his result a
proof.  He writes, ``From personal experience, it is extremely easy to
make programming mistakes. We have taken many precautions, 
\dots
Yet, I want to emphasize that this is only an
experimental result and it desperately needs an independent
verification, or better still, a theoretical
explanation''~\cite{LamS}.

Recent speculation at {\it Math Overflow} holds that the next case,
$n=12$, remains solidly out of computational reach~\cite{Horn}.

\subsection{hyperbolic manifolds}

Computers have helped to resolve a number of open conjectures about
hyperbolic manifolds (defined as complete Riemannian manifolds with
constant negative sectional curvature $-1$), including the proof that
the space of hyperbolic metrics on a closed hyperbolic $3$-manifold is
contractible~\cite{GMT},~\cite{GabICM}.

\subsection{chaos theory and strange attractors}

The theory of chaos has been one of the great success stories of
twentieth century mathematics and science.   Turing\footnote{For
  early history, see \cite[p.~971]{Wolfram:NKS}. Turing vainly hoped that
  digital computers might be insulated from the effects of chaos.}
expressed the notion of chaos with these words, ``quite small errors
in the initial conditions can have an overwhelming effect at a later
time.  The displacement of a single electron by a billionth of a
centimetre at one moment might make the difference between a man being
killed by an avalanche a year later, or escaping''~\cite{Tu50}.
Later, the metaphor became a butterfly
that stirs up a tornado in Texas by flapping its wings in Brazil.

Thirteen years later, Lorenz encountered chaos as he ran weather
simulations on a Royal McBee LGP-30 computer~\cite{Lo63}.  When he
reran an earlier numerical solution with what he thought to be
identical initial data, he obtained wildly different results.  He
eventually traced the divergent results to a slight discrepancy in
initial conditions caused by rounding in the printout.  The {\it Lorenz
  oscillator} is the simplified form of Lorenz's original ordinary
differential equations.

A set $A$ is {\it attracting} if it has a neighborhood $U$ such
that
\[
A = \bigcap_{t\ge 0} f_t(U),
\]
where $f_t(x)$ is the solution of the dynamical system (in present
case the Lorenz oscillator) at time $t$ with initial condition
$x$. That is, $U$ flows towards the attracting set $A$.  Simulations
have discovered attracting sets with strange properties such as
non-integral Hausdorff dimension and the tendency for a small slab of 
volume to quickly spread throughout the attractor.

%%
% XX Define Strange attractor.
% ``a tiny blob of initial values rapidly smears out over the
% entire attractor''~\cite{WT}.

\begin{figure}[h!]
  \centering
\includegraphics[scale=0.3]{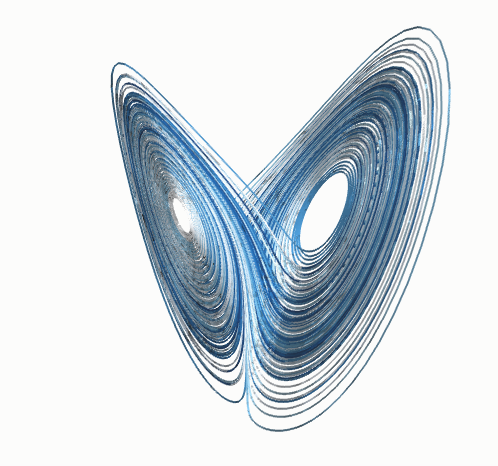}
  \caption{The Lorenz oscillator gives one of the
most famous images of mathematics, a {\it strange attractor} in
dynamical systems (image source~\cite{Lor11}).}
\label{fig:lorenz}
\end{figure}
%\includegraphics[scale=0.7]{200px-Lorenz_attractor_boxed.png}
% Wiki source
% PERM. Wiki commons
% http://en.wikipedia.org/wiki/Lorenz_attractor. 

Lorenz conjectured in 1963 that his oscillator has a strange attractor
(Figure~\ref{fig:lorenz}).  In 1982, the Lax report cited soliton
theory and strange attractors as two prime examples of the ``discovery
of new phenomena through numerical experimentation,'' and calls such
discovery perhaps the most ``significant application of scientific
computing''~\cite{Lax}.  Smale, in his list of $18$ ``Mathematical
Problems for the Next Century'' made the fourteenth problem to present
a rigorous proof that the dynamics of the Lorenz oscillator is a
strange attractor~\cite{Sma98} with various additional properties that
make it a ``geometric Lorenz attractor.''

Tucker has solved Smale's fourteenth problem by
computer~\cite{Tuc02}~\cite{St00}.  One particularly noteworthy aspect
of this work is that chaotic systems, by their very nature, pose
particular hardships for rigorous computer analysis.  Nevertheless,
Tucker implemented the classical Euler method for solving ordinary
differential equations with particular care,  using interval
arithmetic to give mathematically rigorous error bounds.  Tucker has
been awarded numerous prizes for this work, including the Moore Prize
(2002) and the EMS Prize (2004).

Smale's list in general envisions a coming century in which computer
science, especially computational complexity, plays a much larger
role than during the past century. He finishes the list with the
open-ended philosophical problem  that echoes Turing: {\it
  ``What are the limits of intelligence, both artificial and human?''}

\subsection{$4/3$}
\begin{figure}[h!]
  \centering
\includegraphics[scale=0.5]{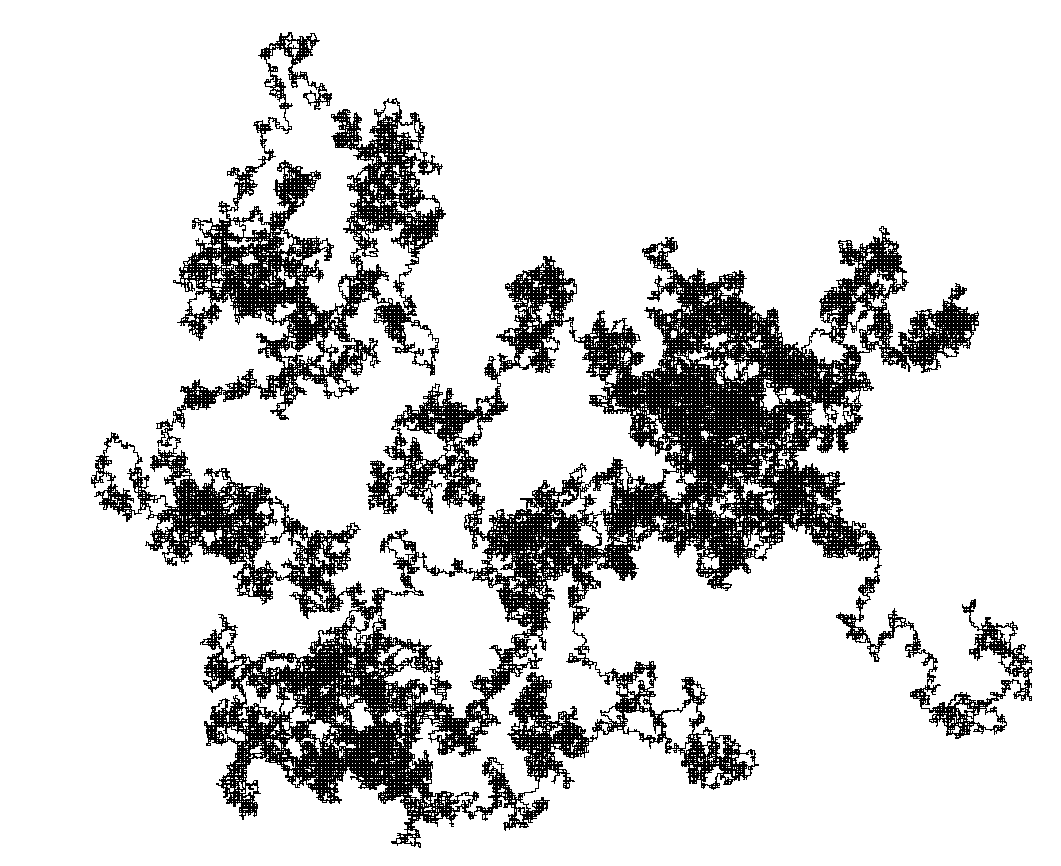}
  \caption{A simulation of planar Brownian motion. Mandelbrot used ``visual
inspection supported by computer experiments'' to formulate deep
conjectures in fractal geometry (image generated from source code at~\cite{LSW01}).}
\label{fig:4/3}
\end{figure}

Mandelbrot's conjectures in fractal geometry have resulted in two
Fields Medals.  Here he describes the discovery of the
$4/3$-conjecture made in~\cite{ManFN}.  ``The notion that these
conjectures might have been reached by pure thought -- with no picture
-- is simply inconceivable.\dots I had my programmer draw a very big
sample [Brownian] motion and proceeded to play with it.''  He goes on
to describe computer experiments that led him to enclose the Brownian
motion into black clusters that looked to him like islands with jagged
coastlines (Figure~\ref{fig:4/3}).
``[I]nstantly, my long previous experience with the coastlines of
  actual islands on Earth came handy and made me suspect that the
  boundary of Brownian motion has a fractal dimension equal to $4/3$''
\cite{Man}.

This conjecture, which  Mandelbrot's trained eye spotted in an instant,
took 18 years to prove~\cite{LSW01}.

%%
% graphic source http://www.mrlonline.org/mrl/2001-008-004/2001-008-004-001.pdf
% No permission needed since original image generated for this publication.

%  ``I was not trying to implement any
%  preconceived idea, simply actively ``fishing'' for new things.\dots
%  In order to achieve homogeneity, I decided to make the motion end
%  where it had started. The resulting motion biting its own tail
%  created a distinctive new shape I call Brownian cluster.  Next the
%  same purely aesthetic consideration led to further processing. The
%  continuing wish to eliminate extraneous complexity made me combine
%  all the points that cannot be reached from infinity without crossing
%  the Brownian cluster. Painting them in black sufficed, once again,
%  to create something quite new, resembling an island.  Instantly, it
%  became apparent that its boundary deserved to be investigated''~\cite{Man}.

%  Mandelbrot uses the phrase ``visual inspection supported by
%  computer experiments'' in
%  http://users.math.yale.edu/users/mandelbrot/web_pdfs/
%  EncycScienceAndTechnologyFractals.pdf
%
%  Mandelbrot's 4/3 conjecture (1982).

\subsection{sphere eversion visualization}

%% Other references:
% http://en.wikipedia.org/wiki/Smale's_paradox
% http://torus.math.uiuc.edu/jms/Papers/isama/eversions.pdf Sullivan
% Levy http://www.math.sunysb.edu/CDproject/OvUM/cool-links/
% www.geom.umn.edu/docs/outreach/oi/history.html

Smale (1958) proved that it is possible to turn a sphere inside out
without introducing any creases.\footnote{I am fond of this example,
because The Scientific American
  article \cite{Phi66} about this theorem was my first exposure to ``real
  mathematics'' as a child.}  For a long time, this paradoxical result
defied the intuition of experts.  R. Bott, who had been Smale's
graduate advisor, refused to believe it at first.  Levy writes that
trying to visualize Smale's mathematical argument ``is akin to
describing what happens to the ingredients of a souffl\'e in minute
detail, down to the molecular chemistry, and expecting someone who has
never seen a souffl\'e to follow this `recipe' in preparing the
dish''~\cite{Le95}.

\begin{figure}[h!]
  \centering
\includegraphics[scale=0.5]{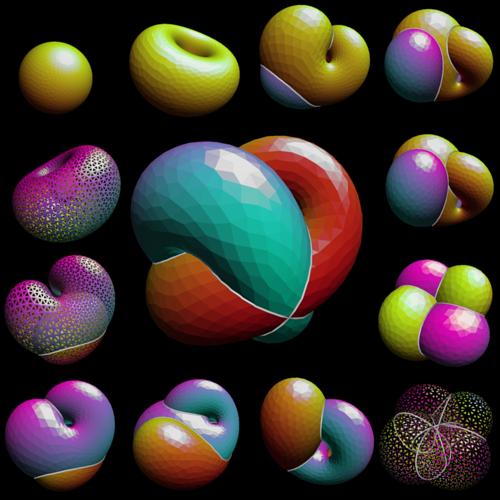}
\caption{Computer-generated stages of a sphere eversion (image
  source~\cite{opti}).}
\label{fig:eversion}
\end{figure}
% http://new.math.uiuc.edu/optiverse/images.html graphics source,
% PERM. 2/3 authors have signed.

It is better to see and taste a souffl\'e first.  The computer videos
of this theorem are spectacular.  Watch them on YouTube!  As we watch
the sphere turn inside out, our intuition grows.  The computer
calculations behind the animations of the first video (the Optiverse)
start with a sphere, half inverted and half right-side out~\cite{SFL}.
From halfway position, the path of steepest descent of an energy
functional is used to calculate the unfolding in both directions to
the round spheres, with one fully inverted
(Figure~\ref{fig:eversion}).  The second video is based on Thurston's
``corrugations''~\cite{LMM}. As the name suggests, this sphere
eversion has undulating ruffles that dance like a jellyfish, but
avoids sharp creases.  Through computers, understanding.

Speaking of Thurston, he contrasts ``our amazingly rich abilities to
absorb geometric information and the weakness of our innate abilities
to convey spatial ideas.\ldots We effortlessly look at a
two-dimensional picture and reconstruct a three-dimensional scene, but
we can hardly draw them accurately''~\cite{Pi11}.  As more and more
mathematics migrates to the computer, there is a danger that
geometrical intuition becomes buried under a logical symbolism.

\subsection{minimal surface visualization} 

\begin{figure}[h!]
  \centering
\includegraphics[scale=0.28]{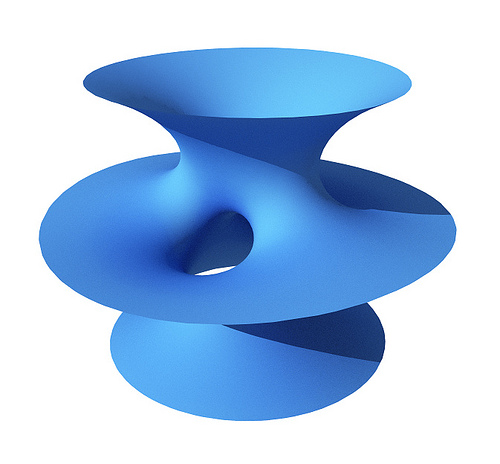}
  \caption{The Costa surface launched an era of computer exploration
in minimal surface theory (image source~\cite{costa}).}
\label{fig:costa}
\end{figure}

Weber and Wolf \cite{WW11} report that the use of computer
visualization has become ``commonplace'' in minimal surface research,
``a conversation between visual aspects of the minimal surfaces and
advancing theory, each supporting the other.''  This started when
computer illustrations of the Costa surface (a particular minimal
surface, Figure~\ref{fig:costa}) in the 1980s revealed dihedral
symmetries of the surface that were not seen directly from its
defining equations.  The observation of symmetry turned out to be the
key to the proof that the Costa surface is an embedding.  The
symmetries further led to a conjecture and then proof of the existence
of other minimal surfaces of higher genus with similar dihedral
symmetries.  As Hoffman wrote about his discoveries, ``The images
produced along the way were the objects that we used to make
discoveries. They are an integral part of the process of doing
mathematics, not just a way to convey a discovery made without their
use'' \cite{Hoffman}.

% http://www.flickr.com/photos/arenamontanus/8059864268/
% PERM:
% Creative commons, attributions & non-commercial license.

\subsection{double bubble conjecture}

Closely related to minimal surfaces are surfaces of constant mean
curvature.  The mean curvature of a minimal surface is zero; surfaces
whose mean curvature is constant are a slight generalization.  They
arise as surfaces that are minimal subject to the constraint that
they enclose a region of fixed volume.  Soap bubble films are
surfaces of constant mean curvature.  

\begin{figure}[h!]
  \centering
\includegraphics[scale=0.15]{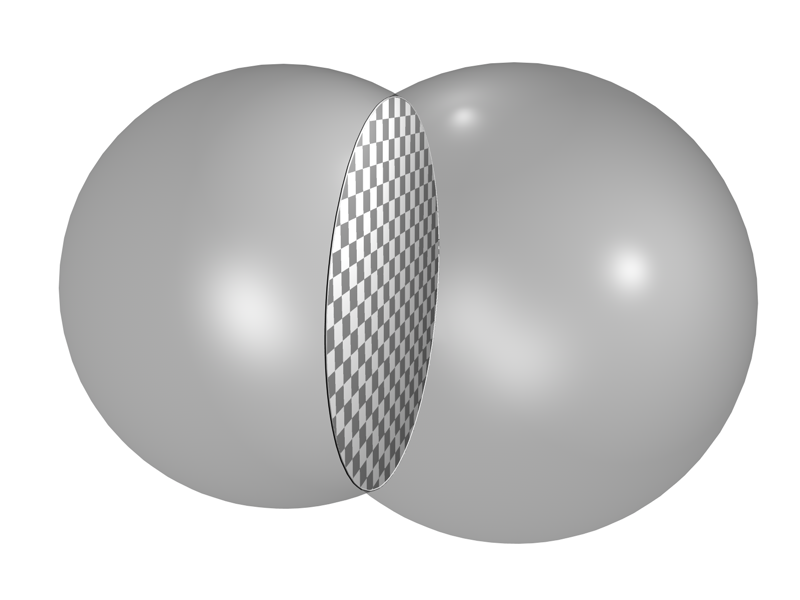}
\caption{The optimality of a double bubble was first established by
  computer, using interval analysis (image source~\cite{double}).}
\label{fig:double}
\end{figure}

The isoperimetric inequality asserts that the sphere minimizes the
surface area among all surfaces that enclose a region of fixed volume.
The double bubble problem is the generalization of the isoperimetric
inequality to two enclosed volumes.  What is the surface minimizing
way to enclose two separate regions of fixed volume?  In the
nineteenth century, Boys~\cite{Boy1890} and Plateau observed
experimentally that the answer should be two partial spherical bubbles
joined along a shared flat disk (Figure~\ref{fig:double}).  The size
of the shared disk is determined by the condition that angles should
be $120^\circ$ where the three surfaces meet.  This is the {\it double
  bubble conjecture}.

The first case of double bubble conjecture to be established was that
of two equal volumes~\cite{HHS95}.  The proof was a combination of
conventional analysis and computer proof.  Conventional analysis
(geometric measure theory) established the existence of a minimizer
and reduced the possibilities to a small number of figures of
revolution, and computers were to analyze each of the cases, showing
in each case by interval analysis either that the case was not a local
minimizer or that its area was strictly larger than the double bubble.
Later theorems proved the double bubble conjecture in the general
unequal volume case without the use of computers~\cite{HMRR}.

% Wiki, Commons.

% old.
% http://math.berkeley.edu/~hutching/pub/bubbles.html 

\begin{figure}[h!]
  \centering
\includegraphics[scale=0.28]{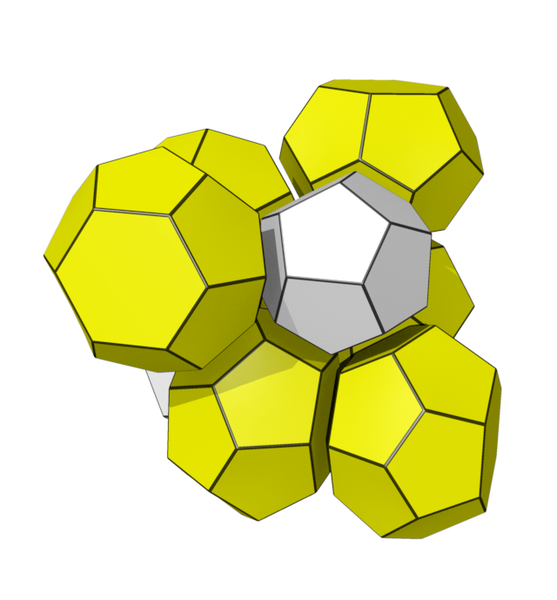}
\caption{The Phelan-Weaire foam, giving the best known partition of
  Euclidean space into cells of equal volume, was constructed with
  Surface Evolver software.  This foam inspired the bubble design of
  the Water Cube building in the 2008 Beijing Olympics (image
  source~\cite{phelan-graphic}).  }
\label{fig:PW}
\end{figure}

The natural extension of the double bubble conjecture from two bubbles
to an infinite bubbly foam is the Kelvin problem.  The problem asks
for the surface area minimizing partition of Euclidean space into
cells of equal volume.  Kelvin's conjecture -- a tiling by slight
perturbations of truncated octahedra -- remained the best known
partition until a counterexample was constructed by two physicists,
Phelan and Weaire in 1993 (Figure~\ref{fig:PW}).  The counterexample
exists not as a physical model, nor as an exact mathematical formula,
but only as an image generated from a triangular mesh in the {\it
  Surface Evolver} computer program.  By default, the counterexample
has become the new conjectural answer to the Kelvin problem, which I
fully expect to be proved someday by computer.

% http://it.wikipedia.org/wiki/File:Foam_-_Weaire-Phelan_structure.png
% wiki graphic, 
% public domain.

\subsection{kissing numbers}

\begin{figure}[h!]
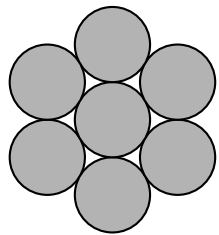

  \centering
\figKISSING
\caption{In two dimensions, the kissing number is $6$. In eight
  dimensions, the answer is $240$.  The proof certificate was found by
  linear programming.}
\label{fig:kissing}
\end{figure}
% generated by Hales.

In the plane, at most six pennies can be arranged in a hexagon so that
they all touch one more penny placed at the center of the hexagon
(Figure~\ref{fig:kissing}).  Odlyzko and Sloane, solved the
corresponding problem in dimension $8$: at most $240$ nonoverlapping
congruent balls can be arranged so that they all touch one more at the
center.

Up to rotation, a unique arrangement of $240$ exists.  To the
cognoscenti, the proof of this fact is expressed as one-line
certificate:
\[
(t - \frac{1}{2})t^2(t + \frac{1}{2})^2 (t + 1).
\]
% polynomial from Pfender-Ziegler.
(For an explanation of the certificates, see \cite{PZ}.)  The
certificate was produced by a linear programming computer search, but
once the certificate is in hand, the proof is computer-free.

As explained above, six is the {\it kissing number} in two dimensions,
$240$ is the kissing number in eight dimensions. In three dimensions,
the kissing number is $12$.  This three-dimensional problem goes back
to a discussion between Newton and Gregory in 1694, but was not
settled until the 1950s.  A recent computer proof makes an exhaustive
search through nearly 100 million combinatorial possibilities to
determine exactly how much the twelve spheres must shrink to
accommodate a thirteenth~\cite{Musin-Tarasov}.  Bachoc and Vallentin
were recently awarded the SIAG/Optimization prize for their use of
semi-definite programming algorithms to establish new proofs of the
kissing number in dimensions $3,4,8$ and new bounds on the kissing
number in various other dimensions~\cite{BV08}.

\subsection{digression on $E_8$}\label{sec:e8}

It is no coincidence that the calculation of Odlyzko and Sloane works
in dimension $8$.  Wonderful things happen in eight dimensional space
and again in $24$ dimensions.

\begin{figure}[h!]
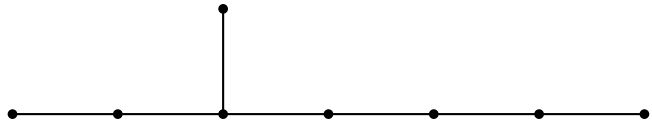

  \centering
\figEE
\caption{The $E_8$ lattice is generated by eight vectors in
  $\ring{R}^8$ whose mutual angles are $120^\circ$ or $90^\circ$
  depending on whether the corresponding dots are joined by a
  segment are not.}
\label{fig:e8}
\end{figure}

Having mentioned the $240$ balls in eight dimensions, I cannot resist
mentioning some further computer proofs.  The centers of the $240$ balls
are vectors whose integral linear combinations generate 
a lattice in $\ring{R}^8$, known as the $E_8$ lattice (Figure~\ref{fig:e8}).

There is a packing of congruent balls in eight dimensions that is
obtained by centering one ball at each vector in the $E_8$ lattice,
making the balls as large as possible without overlap.  Everyone
believes that this packing in eight dimensions is the densest
possible, but this fact currently defies proof.  If the center of the
balls are the points of a lattice, then the packing is called a {\it
  lattice packing}.  Cohn and Kumar have a beautiful computer assisted
proof that the $E_8$ packing is the densest of all lattice packings in
$\ring{R}^8$ (and the corresponding result in dimension $24$ for the
Leech lattice).  The proof is based on the Poisson summation formula.
Pfender and Ziegler's account of this computer-assisted proof won the
Chauvenet Prize of the MAA for writing~\cite{PZ}.

% chauvenet. 2006.
% http://mathdl.maa.org/mathDL/22/?pa=content&sa=viewDocument&nodeId=3065

The $240$ vectors that generate the $E_8$ lattice are the {\it roots}
of a $240+8$ dimensional Lie group (also called $E_8$); that is, a
differentiable manifold that has the analytic structure of a group.
All simple Lie groups were classified in the nineteenth
century.\footnote{I describe the families over $\ring{C}$.  Each
  complex Lie group has a finite number of further real forms.}  They
fall into infinite families named alphabetically, $A_n$, $B_n$, $C_n$,
$D_n$, with $5$ more exceptional cases that do not fall into infinite
families $E_6$, $E_7$ $E_8$, $F_4$, $G_2$.  The exceptional Lie
group\footnote{For decades, $E_8$ has stood for the ultimate in
  speculative physics, whether in heterotic string theory or a
  ``theory of everything.''  Last year, $E_8$ took a turn toward the
  real world, when $E_8$ calculations predicted neutron scattering
  experiments with a cobalt niobate magnet~\cite{BGE8}.} of highest
dimension is $E_8$.

The long-term {\it Atlas Project} aims to use computers to determine
all unitary representations of real reductive Lie groups~\cite{Atlas}.
The $19$-member team focused on $E_8$ first, because everyone respects
the formidable $E_8$.  By 2007, a computer had completed the character
table of $E_8$.  Since there are infinitely many irreducible
characters and each character is an analytic function on (a dense open
subset of) the group, it is not clear without much further explanation
what it might even mean for a computer to output the full character
table as a $60$ gigabyte file~\cite{AdChar}.  What is significant
about this work is that it brings the computer to bear on some
abstract parts of mathematics that have been traditionally largely
beyond the reach of concrete computational description, including
infinite dimensional representations of Lie groups, intersection
cohomology and perverse sheaves.  Vogan's account of this
computational project was awarded the 2011 Conant Prize of the
AMS~\cite{VE8}.

%% Leeuwen's slides describe how a heavily recursive calculation
% requiring 160 GB of RAM was squeezed into 64 GB RAN.
%  (Marc Leeuwen, computational aspects).

% Other references:
% http://atlas.math.umd.edu/software/
% http://atlas.math.umd.edu/talks/cracking-E8.pdf
% http://young.sp2mi.univ-poitiers.fr/~marc/pdf/Zaragoza-E8.pdf
% http://atlas.math.umd.edu/talks/boston.pdf gives 60gb and 77 hours.

While on the topic of computation and representation theory, I cannot
resist a digression into the $P$ versus $NP$ problem, the most
fundamental unsolved problem in mathematics. In my opinion, attempts
to settle $P$ versus $NP$ from the axioms of ZFC are ultimately as
ill-fated as Hilbert's program in the foundations of math (which
nonetheless spurred valuable partial results such as the decision
procedures of Presburger and Tarski), but if I were to place faith
anywhere, it would be in Mulmuley's program in {\it geometric
  complexity theory}.  The program invokes geometric invariant theory
and representation theoretic invariants to tease apart complexity
classes: if the irreducible constituents of modules canonically
associated with two complexity classes are different, then the two
complexity classes are distinct.  In this approach, the determinant
and permanent of a matrix are chosen as the paradigms of what is easy
and hard to compute, opening up complexity theory to a rich
algebro-geometric structure~\cite{Mul11},~\cite{FPNP}.

\subsection{future computer proofs}

Certain problems are natural candidates for computer proof: the Kelvin
problem by the enumeration of the combinatorial topology of possible
counterexamples; the search for a counterexample to the
two-dimensional Jacobian conjecture through the minimal model
program~\cite{Borisov}; resolution of singularities in positive
characteristic through an automated search for numerical quantities
that decrease under suitable blowup; existence of a projective plane
of order $12$ by constraint satisfaction programming; the optimality
proof of the best known packing of tetrahedra in three
dimensions~\cite{Chen-2010}; Steiner's isoperimetric conjecture (1841) for the
icosahedron~\cite{Steiner41}; and the Reinhardt conjecture through
nonlinear optimization~\cite{HR11}.  But proceed with caution!
Checking on our zeal for brute computation, computer-generated
patterns can sometimes fail miserably.  For example, the 
sequence:
\[
\left\lceil{\frac{2}{2^{1/n} - 1}}\right\rceil- 
\left\lfloor{\frac{2 n}{\log 2}}\right\rfloor,\quad n=1,2,3,\ldots
\]
starts out as the zero sequence, but remarkably first gives a nonzero
value when $n$ reaches $777,451,915,729,368$ and then again when
$n=140,894,092,055,857,794$.  See~\cite{Stanley}.
% Stanley example:
% http://mathoverflow.net/questions/15444/the-phenomena-of-eventual-counterexamples
% http://oeis.org/A129935

\smallskip At the close of this first section, we confess that a
survey of mathematics in the age of the Turing machine is a reckless
undertaking, particularly if it almost completely neglects software
products and essential mathematical algorithms -- the Euclidean
algorithm, Newton's method, Gaussian elimination, fast Fourier
transform, simplex algorithm, sorting, Sch\"onhage-Strassen, and many
more.  A starting point for the exploration of mathematical software
is KNOPPIX/Math, a bootable DVD with over a hundred free mathematical
software products (Figure~\ref{fig:knoppix})~\cite{HK08}.  Sage alone
has involved over 200 developers and includes dozens of other
packages, providing an open-source Python scripted alternative to
computer algebra systems such as Maple and Mathematica.

\begin{figure}[h!]
\centering
\begin{tabular}{|@{~~}l@{~~}|@{~~}l@{~~}|}
\hline 
 
\TeX & Active-DVI, AUC\TeX, \TeX{}macs, Kile, Whizzy\TeX\\[0.5ex]
computer algebra & Axiom, CoCoA4, GAP, Macaulay2, Maxima,\\
&PARI/GP, Risa/Asir, Sage, Singular, Yacas\\[0.5ex]
numerical calc&Octave, Scilab, FreeFem++, Yorick\\[0.5ex]
visualization& 3D-XplorMath-J, Dynagraph, GANG, Geomview, \\
 &gnuplot, JavaView, K3DSurf\\[0.5ex]
geometry& C.a.R, Dr.Geo, GeoGebra, GEONExT, KidsCindy, KSEG\\[0.5ex]
programming & CLISP, Eclipse, FASM, Gauche, GCC, Haskell, Lisp\\
&Prolog, Guile, Lazarus, NASM, Objective Caml,\\
&Perl, Python, Ruby, Squeak\\
 [0.5ex]
\hline

\end{tabular}
\caption{Some free mathematical programs on the Knoppix/Math
  DVD~\cite{HK08}.}
\label{fig:knoppix}
\end{figure}

\newpage
\section{Computer Proof}

Proof assistants represent the best effort of logicians, computer
scientists, and mathematicians to obtain complete mathematical rigor
by computer.  This section gives a brief introduction to proof
assistants and describes various recent projects that use them.

The first section described various computer calculations in math, and
this section turns to computer reasoning.  I have never been able to
get used to it being the mathematicians who use computers for
calculation and the computers scientists who use computers for proofs!

\subsection{design of proof assistants}

A formal proof is a proof that has been checked at the level of the
primitive rules of inference all the way back to the fundamental
axioms of mathematics.  The number of primitive inferences is
generally so large that it is quite hopeless to construct a formal
proof by hand of anything but theorems of the most trivial nature.
McCarthy and de Bruijn suggested that we program computers to generate
formal proofs from high-level descriptions of the proof.  This
suggestion has led to the development of proof assistants.

% McC and de Bruijn credited on page 10 of
% ftp://ftp.cs.ru.nl/pub/CSI/CompMath.Found/wiley.pdf

A {\it proof assistant} is an interactive computer program that
enables a user to generate a formal proof from its high-level
description.  Some examples of theorems that have been formally
verified by proof assistants appear in Figure~\ref{fig:table}.  The
computer code that implements a proof assistant lists the fundamental
axioms of mathematics and gives procedures that implement each of the
rules of logical inference.  Within this general framework, there are
enormous variations from one proof assistant to the next. The feature
table in Figure~\ref{fig:table} is reproduced from~\cite{wiedijk:17}.
The columns list different proof assistants, HOL, Mizar, etc.  

Since it is the one that I am most familiar with, my discussion will
focus largely on a particular proof assistant, {\it HOL Light}, which
belongs to the {\it HOL} family of proof assistants. {\it HOL} is an
acronym for Higher-Order Logic, which is the underlying logic of these
proof assistants.  A fascinating account of the history of HOL appears
in \cite{Gor}.  In 1972, R. Milner developed a proof-checking program
based on a deductive system LCF (for Logic of Computable Functions)
that had been designed by Dana Scott a few years earlier.  A long
series of innovations (such as goal-directed proofs and tactics, the ML
language, enforcing proof integrity through the type system,
conversions and theorem continuations, rewriting with discrimination
nets, and higher-order features) have led from LCF to HOL.

\begin{figure}[h!]
\centering
\begin{tabular}{l l l l l}
\hline
Year\hspace{0.5em} &Theorem\hspace{8em} &Proof System\hspace{2em}  &Formalizer\hspace{3em} &Traditional Proof\\ [0.5ex]
\hline \\
1986 &First Incompleteness &Boyer-Moore   &Shankar &G\"odel \\
1990 &Quadratic Reciprocity&Boyer-Moore &Russinoff &Eisenstein\\
1996 &Fundamental - of Calculus &HOL Light &Harrison &Henstock\\
2000 &Fundamental - of Algebra &Mizar &Milewski    &Brynski\\ 
2000 &Fundamental - of Algebra &Coq &Geuvers et al.   &Kneser\\
2004 &Four Color &Coq &Gonthier &Robertson et al.\\
2004 &Prime Number &Isabelle &Avigad et al. &Selberg-Erd\"os\\
2005 &Jordan Curve  &HOL Light &Hales &Thomassen \\
2005 &Brouwer Fixed Point &HOL Light &Harrison &Kuhn \\
2006 &Flyspeck I &Isabelle &Bauer-Nipkow &Hales \\
2007 &Cauchy Residue &HOL Light &Harrison &classical \\
2008 &Prime Number &HOL Light &Harrison &analytic proof \\
2012 &Odd Order Theorem&Coq&Gonthier&Feit-Thompson\\
 [1ex]
\hline
\end{tabular}
\caption{Examples of Formal Proofs, adapted from \cite{Hales:2008:formal}.}
\label{fig:table}
\end{figure}
% Figure from AMS article, Hales, Formal Proof.
% email tchales@gamil.com, Aug 17 from Milewski mentions  Brynski.

Without going into full detail, I will make a few comments about what
some of the features mean.  Different systems can be commended in
different ways: HOL Light for its small trustworthy kernel, Coq for
its powerful type system, Mizar for its extensive libraries, and
Isabelle/HOL for its support and usability.

\bigskip
% typesetting help: 
% http://andrewjpage.com/index.php?/archives/24-Latex-tables-and-rotated-text.html
\noindent
\begin{figure}
\centering
\begin{tabular}{|r|ccccc|ccccc|ccccc|cc|}\hline
~{\it proof assistant}~
&\begin{sideways}HOL\end{sideways}
&\begin{sideways}Mizar\end{sideways}
&\begin{sideways}PVS\end{sideways}
&\begin{sideways}Coq\end{sideways}
&\begin{sideways}Otter/Ivy\end{sideways}
&\begin{sideways}Isabelle/Isar\end{sideways}
&\begin{sideways}Alfa/Agda\end{sideways}
&\begin{sideways}ACL2\end{sideways}
&\begin{sideways}PhoX\end{sideways}
&\begin{sideways}IMPS\end{sideways}
&\begin{sideways}Metamath\end{sideways}
&\begin{sideways}Theorema\end{sideways}
&\begin{sideways}Lego\end{sideways}
&\begin{sideways}Nuprl\end{sideways}
&\begin{sideways}$\Omega$mega\end{sideways}
&\begin{sideways}$B$ method\end{sideways}
&\begin{sideways}Minilog\end{sideways}\\
\hline
small proof kernel (`proof objects')
&+ &- &- &+  &+ &+ &+ &-  &+ &- &+ &-  &+ &- &+ &-  &+ 
\\
calculations can be proved automatically
&+ &- &+ &+  &+ &+ &- &+  &+ &+ &- &+  &+ &+ &+ &+  &+ 
\\
extensible/programmable by the user~
&+ &- &+ &+  &- &+ &- &-  &- &- &- &-  &- &+ &+ &-  &+ 
\\
powerful automation~
&+ &- &+ &-  &+ &+ &- &+  &- &+ &- &+  &- &- &+ &+  &- 
\\
readable proof input files~
&- &+ &- &-  &- &+ &- &+  &- &- &- &+  &- &- &- &-  &- 
\\
\hline
constructive logic supported~
&- &- &- &+  &- &+ &+ &-  &- &- &+ &-  &+ &+ &- &-  &+ 
\\
logical framework~
&- &- &- &-  &- &+ &- &-  &- &- &+ &-  &- &- &- &-  &- 
\\
typed~
&+ &+ &+ &+  &- &+ &+ &-  &+ &+ &- &-  &+ &+ &+ &-  &+ 
\\
decidable types~
&+ &+ &- &+  &- &+ &+ &-  &+ &+ &- &-  &+ &- &+ &-  &+ 
\\
dependent types~
&- &+ &+ &+  &- &- &+ &-  &- &- &- &-  &+ &+ &- &-  &- 
\\
\hline
based on higher order logic~
&+ &- &+ &+  &- &+ &+ &-  &+ &+ &- &+  &+ &+ &+ &-  &- 
\\
based on ZFC set theory~
&- &+ &- &-  &- &+ &- &-  &- &- &+ &-  &- &- &- &+  &- 
\\
large mathematical standard library~
&+ &+ &+ &+  &- &+ &- &-  &- &+ &- &-  &- &+ &- &-  &- 
\\
%statement about $\ring{R}$~
%&+ &+ &+ &+  &- &+ &- &+  &- &+ &+ &+  &- &- &+ &-  &+ 
%\\
%statement about $\sqrt{\phantom X}$~
%&+ &+ &+ &+  &- &+ &- &-  &- &+ &+ &+  &- &- &+ &-  &- 
%\\
\hline
\end{tabular}
\caption{Features of proof assistants~\cite{wiedijk:17}. The table is
published by permission from Springer Science Business Media B.V.}
\label{fig:feature}
\end{figure}
% PERM obtained.
\bigskip 

\subsubsection{small proof kernel.} If a proof assistant is used to
check the correctness of proofs, who checks the correctness of the
proof assistant itself?  De Bruijn proposed that the proofs of a
proof assistant should be capable of being checked by a short piece of
computer code -- something short enough to be checked by hand.  For
example, the kernel of the proof assistant HOL Light is just $430$
lines of very readable computer code.  The architecture of the system
is such that if these $430$ lines are bug free then it is
incapable\footnote{I exaggerate. Section~\ref{sec:trust} goes into
  detail about trust in computers.} of generating a theorem that
hasn't been properly proved.

\subsubsection{automating calculations.} Mathematical argument
involves both calculation and proof.  The foundations of logic often
specify in detail what constitutes a mathematical proof (a sequence of
logical inferences from the axioms), but downgrade calculation to
second-class status, requiring every single calculation to undergo a
cumbersome translation into logic.
Some proof assistants allow {\it reflection} 
(sometimes implausibly attributed to {\it Poincar\'e}), 
which admits as proof the output from a verified
algorithm (bypassing the expansive translation into logic of each
separate execution of the algorithm)~\cite[p.~4]{HPSH},~\cite{BFM}.

%%
% Deleted material:
% Poincar\'e has drawn a distinction between verification of
% individual facts such as $2+2=4$ (which he says ``leads to
% nothing'') and what he calls ``real proofs''~\cite[p.~4]{HPSH}.
% This distinction has come to mean the distinction between the
% algorithmically solvable and the rest. [...] the {\it Poincar\'e
%   principle} holds for that algorithm.  A closely related form of
% {\it reflection} is proof based on the syntactic form of a term,
% rather than the term itself~\cite{BFM}.

\subsubsection{constructive logic.} The law of excluded middle
$\phi\lor\lnot \phi$ is accepted in classical logic, but rejected in
constructive logic.  A proof assistant may be constructive or
classical.  A box ({\it A Mathematical Gem}) shows how HOL Light
becomes classical through the introduction of an axiom of choice.

\newpage
\bigskip
\noindent

\framebox{\parbox{4.6in}{ \smallskip \centerline{\it A Mathematical
      Gem -- Proving the Excluded Middle} \smallskip The logic of HOL
    Light is intuitionistic until the axiom of choice is introduced
    and classical afterwards.  By a result of Diononescu~\cite{Bee85}, choice and
    extensionality imply the law of excluded middle:
 \[\phi\lor \lnot\phi.\]

The proof is such a gem that I have chosen to include it as the only
complete proof in this survey article.  Consider the two sets of booleans
\begin{align*}
P_1 &= \{ x \mid (x = \text{false}) \lor ((x = \text{true})\land\phi)\}\quad \text{and}\\
P_2  &= \{x \mid (x = \text{true}) \lor ((x = \text{false})\land\phi))\}.
\end{align*}
The sets are evidently nonempty, because $\false\in P_1$ and $\true\in P_2$.
By choice, we may pick  $x_1\in P_1$ and $x_2\in P_2$; and by the definition of $P_1$ and $P_2$:
\[
(x_1=\text{false})\lor (x_1 = \text{true}),\qquad (x_2=\text{false})\lor (x_2 = \text{true}). 
\]
We may break the proof of the excluded middle into four cases, depending on the
two possible truth values of each of $x_1$ and $x_2$.

{\bf Cases} $(x_1,x_2)=(\true,\true),$ $(x_1,x_2)=(\true,\false)$: By
the definition of $P_1$, if $x_1 = \true$, then $\phi$, so $\phi\lor
\lnot\phi$. 

 {\bf Case} $(x_1,x_2)=(\false,\false)$: Similarly, by the definition
of $P_2$, if $x_2=\false$, then $\phi$, so also $\phi\lor\lnot\phi$.

{\bf Case} $(x_1,x_2)=(\false,\true)$: If $\phi$, then $P_1=P_2$, and
the choices $x_1$ and $x_2$ reduce to a single choice $x_1=x_2$, which
contradicts $(x_1,x_2)=(\false,\true)$.  Hence $\phi$ implies
$\false$; which by the definition of negation gives $\lnot\phi$, so
also $\phi\lor\lnot\phi$.\hfill{\it Q.E.D.}

}}

\newpage

\subsubsection{logical framework.} 
Many different systems of logic arise in computer science.  In some
proof assistants the logic is fixed. Other proof assistants are more
flexible, allowing different logics to be plugged in and played with.
The more flexible systems implement a meta-language, a {\it logical
  framework}, that gives support for the implementation of multiple
logics.  Within a logical framework, the logic and axioms of a proof
assistant can themselves be formalized, and machine
translations\footnote{My long term Flyspeck project seeks to give a
  formal proof of the Kepler conjecture~\cite{Hales:2005:DSP}. This
  project is now scattered between different proof assistants. Logical
  framework based translations between proof assistants gives me hope
  that an automated tool may assemble the scattered parts of the
  project.} can be constructed between different foundations of
mathematics~\cite{IRFF}.

\subsubsection{type theory.}

Approaching the subject of formal proofs as a mathematician whose
practice was shaped by Zermelo-Fraenkel set theory, I first treated
types as nothing more than convenient identifying labels (such real
number, natural number, list of integers, or boolean) attached to
terms, like the
PLU stickers on fruit that get peeled away before consumption.  Types
are familiar from programming languages as a way of identifying what
data structure is what.  In the simple type system of HOL Light, to
each term is affixed a unique type, which is either a primitive type
(such as the boolean type $\op{\it bool}$), a type variable
($A,B,C,\ldots$), or inductively constructed from other types with the
arrow constructor ($A\rightarrow B$, $A\rightarrow (\op{\it bool}
\rightarrow C)$, etc.).  There is also a way to create subtypes of
existing types.  If the types are interpreted naively as sets, then
$x\tc A$ asserts that the term $x$ is a member of $A$, and
$f:A\rightarrow B$ asserts that $f$ is a member of $A\rightarrow B$,
the set of functions from $A$ to $B$.

In untyped set theory, it is possible to ask ridiculous questions such
as whether the real number $\pi=3.14\ldots$, when viewed as a raw set,
is a finite group.  In fact, in a random exploration of set theory,
like a monkey composing sonnets at the keyboard, ridiculous questions
completely overwhelm all serious content.  Types organize data on the
computer in meaningful ways to cut down on the static noise in the
system.  The question about $\pi$ and groups is not well-typed and
cannot be asked.  Russell's paradox also disappears: $X \not\in X$ is
not well-typed.  For historical reasons, this is not surprising:
Russell and Whitehead first introduced types to overcome the paradoxes of
set theory, and from there, through Church, they passed into computer
science.

Only gradually have I come to appreciate the significance of a
comprehensive {\it theory of types}.  The type system used by a proof
assistant determines to a large degree how much of a proof the user
must contribute and how much the computer automates behind the scenes.
The type system is {\it decidable} if there is a decision procedure to
determine the type of each term.

A type system is {\it dependent} if a type can depend on 
another term.  For example, Euclidean space $\ring{R}^n$, depends on
its dimension $n$.  For this reason, Euclidean space is most naturally
implemented in a proof assistant as a dependent type.  
In a proof assistant such as HOL Light that does not have dependent types,
extra work is required to develop a Euclidean space library.

\subsection{propositions as types}

I mentioned the naive interpretation of each type $A$ as a set and a
term $x\tc A$ as a member of the set.  A quite different
interpretation of types has had considerable influence in the design
of proof assistants.  In this ``terms-as-proofs'' view, a type $A$
represents a proposition and a term $x\tc A$ represents a proof of the
proposition $A$.  A term with an arrow type, $f:A\rightarrow B$, can
be used to to construct a proof $f(x)$ of $B$ from a proof $x$ of $A$.
In this interpretation, the arrow is logical implication.

A further interpretation of types comes from programming languages.
In this ``terms-as-computer-programs'' view, a term is a program
and the type is its specification.   For example,
$f:A\rightarrow B$ is a program $f$ that takes input of type
$A$ and returns a value of type $B$.

By combining the ``terms as proofs'' with the ``terms as computer
programs'' interpretations, we get the famous {\it Curry-Howard
  correspondence} that identifies proofs with computer programs and
identifies each proposition with the type of a computer program. 
For example,  the most fundamental rule of logic,
\[
\frac{A,\quad A\rightarrow B}{B},  \qquad \text{\it(modus ponens)}
\]
(from $A$ and $A\text{-implies-}B$ follows $B$) is identified with the
function application in a computer program; from $x\tc A$ and
$f:A\rightarrow B$ we get $f(x)\tc B$.  To follow the correspondence
is to extract an executable computer program from a mathematical
proof.  The Curry-Howard correspondence has been extremely fruitful,
with a multitude of variations, running through a gamut of proof
systems in logic and identifying each with a suitable programming
domain.

\subsection{proof tactics}

In some proof assistants, the predominant proof style is a backward
style proof.  The user starts with a {\it goal}, which is a statement
to be proved. In interactive steps, the user reduces the goal to
successively simpler goals until there is nothing left to prove.

Each command that reduces a goal to simpler goals is called a {\it
  tactic}.  For example, in the proof assistant HOL Light, there are
about 100 different commands that are tactics or higher-order
operators on tactics (called tacticals).  Figure~\ref{fig:tactic}
shows the most commonly used proof commands in HOL Light.  The most
common tactic is {\it rewriting}, which takes a theorem of the form
$a=b$ and substitutes $b$ for an occurrence of $a$ in the goal.

\bigskip
\noindent
\begin{figure}
\centering
\begin{tabular}{|@{~~}l@{~}|@{~}l@{~}|@{~}c@{~~}|}\hline
\text{\it name}  &~\text{\it purpose} &\text{\it usage}\\
\hline
\text{THEN}   &~\text{combine two tactics into one}   & 37.2\%\\
\text{REWRITE} &~\text{use $a=b$ to replace $a$ with $b$ in goal} & 14.5\%\\
\text{MP\_TAC} &~\text{introduce a previously proved theorem} &4.0\%\\
\text{SIMP\_TAC}&~\text{rewriting with conditionals} & 3.1\%\\
\text{MATCH\_MP\_TAC} &~\text{reduce a goal $b$ to $a$, given a theorem $a\Longrightarrow b$}& 3.0\%\\
\text{STRIP\_TAC} &~\text{(bookkeeping) unpackage a bundled goal} & 2.9\%\\
\text{MESON\_TAC}&~\text{apply first-order reasoning to solve the goal} & 2.6\%\\
\text{REPEAT} &~\text{repeat a tactic as many times as possible} & 2.5\%\\
\text{DISCH\_TAC}&~\text{(bookkeeping) move hypothesis to the assumption list\!\!\!} & 2.3\%\\
\text{EXISTS\_TAC}&~\text{instantiate an existential goal $\exists x\dots$}& 2.3\%\\
\text{GEN\_TAC}&~\text{instantiate a universal goal $\forall x\dots$}& 1.4\%
\\
\hline
\end{tabular}
\caption{A few of the most common proof commands in the HOL Light proof assistant}
\label{fig:tactic}
\end{figure}
\bigskip

In the Coq proof assistant, the tactic system has been streamlined to
an extraordinary degree by the {\it SSReflect} package, becoming a
model of efficiency for other proof assistants to emulate, with an
extremely small number of tactics such as the {\it move} tactic for
bookkeeping, one for rewriting, ones for forward and backward
reasoning, and another for case
analysis~\cite{gonISSR},~\cite{gonSSRE}.  The package also provides
support for exploiting the computational content of proofs, by
integrating logical reasoning with efficient computational algorithms.

\subsection{first-order automated reasoning}

Many proof assistants support some form of automated reasoning to
relieve the user of doing rote logic by hand. For example,
Table~\ref{fig:tactic} lists {\it meson} (an acronym for Loveland's
Model Elimination procedure), which is HOL Light's tactic for
automated reasoning~\cite[Sec.~3.15]{Ha09}, ~\cite{harrison:meson}.
The various automated reasoning tools are generally {\it first-order}
theorem provers.  The classic resolution algorithm for first-order
reasoning is illustrated in a box ({\it Proof by Resolution}).

\newpage
\bigskip
\noindent

\framebox{\parbox{4.6in}{ \smallskip \centerline{\it Proof by
Resolution} \smallskip 
\parskip=\baselineskip

Resolution is the granddaddy of automated reasoning in first-order
logic.  The resolution rule takes two disjunctions
\[P\lor A \qquad  \text{and}\qquad \lnot P' \lor B\] and
concludes \[A'\lor B',\] where $A'$ and $B'$ are the specializations of
$A$ and $B$, respectively, under the {\it most general unifier} of $P$ and $P'$.
(Examples of this in practice appear below.)

This box presents a rather trivial example of proof by resolution, to
deduce the easy theorem asserting that every infinite set has a
member.  The example will use the following notation.  Let $\emptyset$
be a constant representing the empty set and constant $c$ representing
a given infinite set. We use three unary predicates $e$, $f$, $i$ that
have interpretations
\[
e(X)~~\text{``$X$ is empty''},
\quad f(X)~~\text{``$X$ is finite''}, 
\quad i(X)~~\text{``$X$ is infinite.''}
\]
The binary predicate $(\in)$ denotes set membership.
We prove $i(c) \Rightarrow (\exists z. z\in c)$ 
``an infinite set has a member'' by resolution. 

To argue by contradiction,
we introduce the
hypothesis $i(c)$ and the negated conclusion $\lnot (Z\in c)$ as axioms.
Here are the axioms that we allow in the deduction.  The axioms have
been preprocessed, stripped of quantifiers, and written as a disjunction of literals.
Upper case letters are variables.

\smallskip\noindent
\begin{tabular}{lll}
&{\it Axiom}&{\it Informal Description}\\
1.~&$i(c)$~~&Assumption of desired theorem.\\
2.&$\lnot (Z\in c)$~~&Negation of conclusion of desired theorem.\\
3.&$e(X) \lor (u(X)\in X)$~~&A nonempty set has a member.\\
4.&$e(\emptyset)$~~&The empty set is empty.\\
5.&$f(\emptyset)$~~&The empty set is finite.\\
6.&$\lnot i(Y) \lor \lnot f(Y)$~~&A set is not both finite and infinite.\\
7.& $\lnot e(U) \lor \lnot e(V) \lor \lnot i(U) \lor i(V)$~~&Weak indistinguishability of empty sets.\\
\end{tabular}
\smallskip

Here are the resolution inferences from this list of axioms. The final step obtains the
desired contradiction. 

\smallskip
\begin{tabular}{lll}
{\it }&{\it Inference}&{\it Resolvant}\\
~8.&(resolving 2,3, unifying $X$ with $c$ and $u(X)$ with $Z$)~~~~~&$e(c)$\\
~9.&(resolving 7,8, unifying $U$ with $c$)~~~&$\lnot e(V) \lor \lnot i(c) \lor i(V)$\\
10.&(resolving 1,9)~~~&$\lnot e(V) \lor  i(V)$\\
11.&(resolving 4,10, unifying $V$ with $\emptyset$)~~~&$i(\emptyset)$\\
12.&(resolving 6,11, unifying $Y$ with $\emptyset$)~~~&$\lnot f(\emptyset)$\\
13.&(resolving 12,5)~~~&$\perp$\\
\end{tabular}

{{\it Q.E.D.}}

}}

\newpage

Writing about first-order automated reasoning, Huet and Paulin-Mohring
\cite{Coq} describe the situation in the early 1970s as a
``catastrophic state of the art.''  ``The standard mode of use was to
enter a conjecture and wait for the computer's memory to exhaust its
capacity.  Answers were obtained only in exceptionally trivial
cases.'' 
They go on to describe numerous developments (Knuth-Bendix, LISP,
rewriting technologies, LCF, ML, Martin-L\"of type theory, NuPrl,
Curry-Howard correspondence, dependent types, etc.) that led up to the
Coq proof assistant.  These developments led away from first-order
theorem proving with its ``thousands of unreadable logical
consequences'' to a highly structured approach to theorem proving in
Coq.

%%
% translation by T. Hales from ``Le mode standard d'utilisation
% \'etait de rentrer sa conjecture et d'attendre que la m\'emoire de
% l'ordinateur soit satur\'ee.  Seulement dans des cas
% exceptionnellement triviaux une r\'eponse \'etait obtenue.''

First-order theorem proving has developed significantly over the years
into sophisticated software products.  They are no longer limited to
``exceptionally limited cases.''   Many different
software products compete in an annual competition (CASC), to see
which can solve difficult first-order problems the fastest.  The LTB
(large theory batch) division of the competition includes problems
with thousands of axioms~\cite{PSST}.  Significantly, this is the same order of
magnitude as the total number of theorems in a proof assistant.  What
this means is that a first-order theorem provers have reached the
stage of development that they might be able to give fully automated
proofs of new theorems in a proof assistant, working from the full
library of previously proved theorems.

\subsubsection{sledgehammer.}

The Sledgehammer tactic is Paulson's implementation of this idea of full automation in
the Isabelle/HOL proof assistant~\cite{Paar}.  As the name `Sledgehammer' suggests,
the tactic is all-purpose and powerful, but demolishes all higher
mathematical structure, treating every goal as a massive unstructured
problem in first-order logic.  If $L$ is the set of all theorems in
the Isabelle/HOL library, and $g$ is a goal, it would be possible to hand
off the problem $L\Longrightarrow g$ to a first-order theorem
prover.  However, success rates are dramatically improved, when the
theorems in $L$ are first assessed by heuristic rules for their likely
relevance for the goal $g$, in a process called {\it relevance
  filtering}. 
This filtering is used to reduce $L$ to an axiom set $L'$ of a
few hundred theorems that are deemed most likely to prove $g$.

The problem $L'\Longrightarrow g$ is stripped of type information,
converted to a first-order, and fed to
first-order theorem provers. Experiments indicate that
it is more effective to feed a problem in parallel into multiple
first-order provers for a five-second burst than to hand the problem
to the best prover (Vampire) for a prolonged attack~\cite{Paar},~\cite{Boehme-Nipkow-IJCAR10}.
%
% ``B\"ohme and Nipkow \cite{Boehme-Nipkow-IJCAR10} have demonstrated that
%running three different [first-order] theorem provers (E, SPASS and
%Vampire) for five seconds solves as many problems as running the
%best theorem prover (Vampire) for two full minutes.  It would be
%better to utilise even more theorem provers'' 
When luck runs in your favor, one of the first-order theorem provers
finds a proof.

The reconstruction of a formal proof from a first-order proof can
encounter hurdles.  For one thing, when type information is stripped
from the problem (which is done to improve performance), soundness is
lost.  ``In unpublished work by Urban, MaLARea [a machine learning
program for relevance ranking] easily proved the full Sledgehammer
test suite by identifying an inconsistency in the translated lemma
library; once MaLARea had found the inconsistency in one proof, it
easily found it in all the others'' \cite{Paar},~\cite{UrM}.
Good results have been obtained in calling the first-order prover
repeatedly to find a smaller set of axioms $L''\subset L'$ that imply
the goal $g$. A manageably sized set $L''$ is then passed to the metis
tactic\footnote{Metis is a program that automates first-order
  reasoning~\cite{Metis}.} in Isabelle/HOL, which
constructs a formal proof $L''\Longrightarrow g$ from scratch.

% Metis reference:
% http://www.gilith.com/software/metis/notes.html

B\"ohme and Nipkow 
took 1240 proof goals that appear in several diverse theories of the
Isabelle/HOL system and ran sledgehammer on all of
them~\cite{Boehme-Nipkow-IJCAR10}. The results are astounding. The
success rate (of obtaining fully reconstructed formal proofs) when
three different first-order provers run for two-minutes each was 48\%.
The proofs of these same goals by hand might represent years of human
labor, now fully automated through a single new tool.

Sledgehammer has led to a new style of theorem proving, in which the
user is primarily responsible for stating the goals.  In the final
proof script, there is no explicit mention of sledgehammer.  Metis
proves the goals, with sledgehammer operating silently in the
background to feed metis with whatever theorems it needs.  For
example, a typical proof script might contain lines such as
\cite{Paar}
\[
\begin{array}{ll}
\text{{\bf hence} ``$x \subseteq \text{space}~M$''}\\
\text{~~~{\bf by} (metis sets into space lambda system sets)}&
\end{array}
\]
The first line is the goal that the user types. The second line has been
automatically inserted into the proof script by the system, with the relevant
theorems {\tt sets, into} etc. selected by Sledgehammer.

\subsection{computation in proof assistants.}

One annoyance of formal proof systems is the difficulty in locating
the relevant theorems.  At last count, HOL Light had about $14,000$
theorems and nearly a thousand procedures for proof construction.
Larger developments, such as Mizar, have about twice as many theorems.
Good search tools have somewhat relieved the burden of locating
theorems in the libraries.  However, as the formal proof systems
continue to grow, it becomes ever more important to find ways to use
theorems without mentioning them by name.

% HOL Light has 933 indexed help entries and
% 14,030 indexed theorems as of June 12, 2011 with flyspeck loaded.

As an example of a feature which commendably reduces the burden of
memorizing long lists of theorem names, I mention the {REAL\_RING}
command in HOL Light, which is capable of proving any system of
equalities and inequalities that holds over an arbitrary integral
domain.  For example, I can give a one-line formal proof of  an isogeny
$(x_1,y_1) \mapsto (x_2,y_2)$ of elliptic curves: if we have a point
on the first elliptic curve:
\begin{align*}
y_1^2 &= 1 + a x_1^2 + b x_1^4,\\
x_2 y_1&=x_1,\\
y_2 y_1^2&=(1 - b x_1^4),\\  %&1 - b * x1 pow 4 
y_1&\ne 0
\end{align*}
then $(x_2,y_2)$ lies on a second elliptic curve
\[
y_2^2 = 1 + a' x_2^2 + b' x_2^4,  
\]
where $a' = -2a$ and $b' = a^2 - 4b$.  In the proof assistant, 
the input of the
statement is as economical as what I have written here. We expect computer
algebra systems to be capable of checking identities like this, but to
my amazement, I found it {\it easier} to check this isogeny in HOL
Light than to check it in {\it Mathematica}.

The algorithm works in the following manner.  A universally quantified
system of equalities and inequalities holds over all integral domains
if and only if it holds over all fields.  By putting the formula in
conjunctive normal form, it is enough to prove a finite number of
polynomial identities of the form:
\begin{equation}\label{eqn:q}
(p_1=0) \lor \cdots \lor (p_n=0) \lor (q_1\ne 0) \lor\cdots\lor (q_k\ne 0).
\end{equation}
An element in a field is zero, if and only if it is not a unit.  Thus
we may rewrite each polynomial equality $p_i=0$ as an equivalent
inequality $1-p_i z_i\ne 0$.  Thus, without loss of generality, we may
assume that $n=0$; so that all disjuncts are inequalities.  The
formula (\ref{eqn:q}) is logically equivalent to
\[
(q_1 =0) \land \cdots \land (q_k = 0) \Longrightarrow \text{false}.
\]
In other words, it is enough to prove that the zero set of the ideal
$I=(q_1,\ldots,q_n)$ is empty.  For this, we may use
Gr\"obner bases\footnote{Kaliszyk's benchmarks suggest that the
  Gr\"obner basis algorithm in the proof assistant Isabelle runs about
  twenty times faster than that of HOL Light.
}  to prove that $1\in I$, to certify that
the zero set is empty.  
%% reference: CICM 2011 presentation by K.

%% REAL_RING `a' = &2 * a /\ b' = a*a - &4 * b /\ 
% x2 * y1 = x1 /\ y2 * y1 pow 2 = &1 - b * x1 pow 4 /\ 
% y1 pow 2 = &1 + a * x1 pow 2 + b * x1 pow 4 /\ ~(y1 = &0) ==> 
% y2 pow 2 = &1 - a' * x2 pow 2 + b' * x2 pow 4`;;

Gr\"obner basis algorithms give an example of a {\it
  certificate-producing procedure}.  A formal proof is obtained in two
stages.  In the first stage an unverified algorithm produces a
certificate.  In the second stage the proof assistant analyzes the
certificate to confirm the results.  Certificate-producing procedures
open the door to external tools, which tremendously augment the power
of the proof assistant.  The meson is procedure implemented this way,
as a search followed by verification.  Other certificate-producing
procedures in use in proof assistants are linear programming, SAT, and
SMT.

\smallskip

Another praiseworthy project is Kaliszyk and Wiedijk's implementation
of a computer algebra system on top of the proof assistant HOL Light.
It combines the ease of use of computer algebra with the rigor of
formal proof~\cite{kaliszyk_p04_calc}.  Even with its notational
idiosyncrasies (\verb!&! and \verb!#! as a markers of real numbers,
\verb!Cx! as a marker of complex numbers, \verb!ii! for $\sqrt{-1}$,
and \verb!--! for unary negation), it is the kind of product that I
can imagine finding widespread adoption by mathematicians. Some of the features
of the system are shown in Figure~\ref{fig:kw}.

\begin{figure}
\begin{verbatim}
 In1 := (3 + 4 DIV 2) EXP 3 * 5 MOD 3 
Out1 := 250 
 In2 := vector [&2; &2] - vector [&1; &0] + vec 1 
Out2 := vector [&2; &3] 
 In3 := diff (diff (\x. &3 * sin (&2 * x) + &7 + exp (exp x))) 
Out3 := \x. exp x pow 2 * exp (exp x) + exp x * exp (exp x) + -- &12 * sin (&2 * x) 
 In4 := N (exp (&1)) 10 
Out4 := #2.7182818284 + ... (exp (&1)) 10 F 
 In5 := 3 divides 6 /\ EVEN 12 
Out5 := T 
 In6 := Re ((Cx (&3) + Cx (&2) * ii) / (Cx (-- &2) + Cx (&7) * ii)) 
Out6 := &8 / &53 
\end{verbatim}
\caption{Interaction with a formally verified computer algebra system~\cite{kaliszyk_p04_calc}.}
\label{fig:kw}
\end{figure}

\subsection{formalization of finite group theory}

The Feit-Thompson theorem, or odd-order theorem, is one of the most
significant theorems of the twentieth century.  (For his work,
Thompson was awarded the three highest honors in
the mathematical world: the Fields Medal, the Abel Prize, and the Wolf
Prize.)  The Feit-Thompson theorem states that every finite simple
group has even order, except for cyclic groups of prime order.  The
proof, which runs about 250 pages, is extremely technical.  The
Feit-Thompson theorem launched the endeavor to classify all finite
simple groups, a monumental undertaking that consumed an entire
generation of group theorists.

Gonthier's team has formalized the proof of the Feit-Thompson theorem~\cite{gon-FT}.
To me as a mathematician, nothing else that has been done by the formal
proof community compares in splendor to the formalization of this
theorem.  Finally, we are doing real mathematics!  The project
formalized two books, \cite{BG94} and \cite{P00}, as well as a significant
body of background material.

The structures of abstract algebra -- groups, rings, modules,
algebras, algebraically closed fields and so forth -- have all been
laid out formally in the Coq proof assistant.  Analogous algebraic
hierarchies appear in systems such as OpenAxiom, MathScheme, Mizar, and
Isabelle; and while some of these hierarchies are elaborate, none have
delved so deeply as the development for Feit-Thompson.
It gets multiple abstract structures to work coherently
together in a formal setting. ``The problem is not so much in capturing
the semantics of each individual construct but rather in having all
the concepts working together well''~\cite{gonMF}. 

% references:
% http://www.open-axiom.org/,
% http://www.cas.mcmaster.ca/research/mathscheme/ (Carette)

\begin{figure}
{

\obeylines\tt
Structure~finGroupType~Type~:= FinGroupType \{
~~~element~:>~finType;
~~~~~~~1~:~element;
~~~~~~\hskip0.8mm ${}^{-1}$~:~element $\to$ element;
~~~~~~~*~:~element $\to$ element $\to$ element;
~~~unitP~:~$\forall\,x,~1*x = x$;
~~~~invP~:~$\forall\,x,~x^{-1} * x = 1$;
~~~~mulP~:~$\forall\,x_1~x_2~x_3,~ x_1 * (x_2 * x_3) = (x_1 * x_2) * x_3$
\}.

}
\caption{The structure of a finite group~\cite{gonMF}.}
\label{fig:group}
\end{figure}

The definition of a finite group in Coq is similar to the textbook definition,
expressed in types and structures (Figure~\ref{fig:group}).
It declares a finite type called {\tt element} that is the group
carrier or domain.  The rest of the structure specifies a left-unit
element $1$, a left-inverse ${}^{-1}$ and an associative binary
operation $( * )$.  

Other aspects of Gonthier's recent work can be found at
\cite{gonPFSF}, \cite{gonPMS}, \cite{gonC}.
Along different lines, a particularly elegant organization of abstract
algebra and category theory is obtained with type
classes~\cite{SpE11}.

\subsection{homotopy type theory}

The simple type theory of HOL Light is adequate for real analysis,
where relatively few types are needed -- one can go quite far with
natural numbers, real numbers, booleans, functions between these
types, and a few functionals.  However, the dependent type theory of
Coq is better equipped than HOL Light for the hierarchy of structures
from groups to rings of abstract algebra.  But even
Coq's type theory is showing signs of strain in dealing with abstract
algebra.  For instance, an unpleasant limitation of Coq's theory of
types is that it lacks the theorem of extensionality for functions: if
two functions take the same value for every argument, it {\it does
  not} follow that the two functions are equal.\footnote{HOL Light
  avoids this problem by positing extensionality as a mathematical axiom.}  The
gymnastics to solve the problem of function extensionality in the
context of the Feit-Thompson theorem are found in ~\cite{gonMF}.

A lack of function extensionality is an indication
that equality in type theory may be misconceived. Recently, {\it
  homotopy type theory} has exploded onto the scene, which turns to
homotopy theory and higher categories as models of type
theory~\cite{htt}.  It is quite natural to interpret a dependent type
(viewed as a family of types parametrized by a second type)
topologically as a fibration (viewed as a family of fibers
parametrized by a base space)~\cite{AW09}. Voevodsky took the
homotopical notions of equality and equivalence and translated them
back into type theory, obtaining the {\it univalence axiom} of type
theory, which posits what types are equivalent~\cite{VV11}, \cite{PW12}, 
\cite{KLV12}, \cite{KLV12a}.  One
consequence of the univalence axiom is the theorem of extensionality
for functions.  Another promising sign for computer theorem-proving
applications is that the univalence axiom appears to preserve the
computable aspects of type theory (unlike for instance, the axiom of
choice which makes non-computable choices)~\cite{LH11}.  We may hope
that some day there may be a back infusion of type-theoretic proofs
into homotopy theory.

\subsection{language of mathematics}

Ganesalingam's thesis is the most significant linguistic study of the
language of mathematics to date~\cite{Gan09}, \cite{Gan10}.
Ganesalingam was awarded the 2011 Beth Prize for the best dissertation
in Logic, Language, or Information.  Although this research is still
at an early stage, it suggests that the mechanical translation of
mathematical prose into formal computer syntax that faithfully
represents the semantics is a realistic hope for the not-to-distant
future.  

% references:
% http://www.cl.cam.ac.uk/news/2011/07/mohan-ganesalingam-awarded-beth-prize/

The linguistic problems surrounding the language of mathematics differ
in various ways from those of say standard English.  A mathematical
text introduces new definitions and notations as it progresses,
whereas in English, the meaning of words is generally fixed from the
outset.  Mathematical writing freely mixes English with symbolic
expressions.  At the same time, mathematics is self-contained in a way
that English can never be; to understand English is to understand the
world.  By contrast, the meaning in a carefully written mathematical
text is determined by Zermelo-Fraenkel set theory (or your favorite
foundational system).

Ganesalingam's analysis of notational syntax is general enough to
treat quite general mixfix operations generalizing infix (e.g. +),
postfix (e.g. factorial !), and prefix ($\cos$).  He analyzes
subscripted infix operators (such as a semidirect product
$H\rtimes_\alpha N$), multi-symboled operators (such as the
three-symboled $[~:~]$ operator for the degree $[K:k]$ of a
field-extension), prefixed words ($R$-module), text within formulas
$\{(a,b) \mid a \text{~is a factor of~} b\}$, unusual script placement
${}^LG$, chained relations $a<b<c$, ellipses $1+2+\cdots+n$,
contracted forms $x,y\in\ring{N}$, and exposed formulas (such as ``for
all $x>0$, \dots'' to mean ``for all $x$,~if $x>0$, then \dots'').

The thesis treats what is called the formal mode of the language of
mathematics -- the language divested of all the informal side-remarks.
The syntax is treated as a context-free grammar, and the semantics are
analyzed with a variant of {\it discourse representation theory},
which in my limited understanding is something very similar to
first-order logic; but different in one significant aspect: it
provides a theory of pronoun references; or put more precisely, a
theory of what may be the ``legitimate antecedent for anaphor.''

A major issue in Ganesalingam's thesis is the resolution of ambiguity.
For example, in the statement
\begin{equation}\label{eqn:P}
P \text{\it ~is prime}
\end{equation}
the term `prime' may mean prime number, prime ideal, or prime
manifold.  His solution is to attach type information to terms (in the
sense of types as discussed above).  The reading of (\ref{eqn:P})
depends on the type of $P$, variously a number, a subset of a ring, or
a manifold.  In this analysis, resolution of ambiguity becomes a task
of a type inference engine.

Because of the need for type information, Ganesalingam raises
questions about the suitability of Zermelo-Fraenkel set theory as the
ultimate semantics of mathematics.  A number of formal-proof
researchers have been arguing in favor of typed foundational systems
for many years.  It is encouraging that there is remarkable
convergence between Ganesalingam's linguistic analysis, innovations in
the Mizar proof assistant, and the development of abstract algebra in
Coq. For example, in various camps we find ellipses (aka big
operators), mixfix operators, type inference, missing argument
inference mechanisms, and so forth.  Also see~\cite{Hoe11}
and~\cite{Forthel}.  Mathematical abuses of notation have turned out
to be rationally construed after all!

\subsection{looking forward}

Let's take the long term view that the longest proofs of the last
century are of insignificant complexity compared to what awaits.  Why
would we limit our creative endeavors to $10,000$ page proofs when we
have tools that allow us to go to a million pages or more?  So far it
is rare for a computer proof has defied human understanding.  No human
has been able to make sense of an unpublished 1500 page
computer-generated proof about Bruck loops\footnote{The theorem states
  that Bruck loops with abelian inner mapping group are centrally
  nilpotent of class two.}~\cite{Stan}.  Eventually, we will have to
content ourselves with fables that approximate the content of a
computer proof in terms that humans can comprehend.

Turing's great theoretical achievements were to delineate what a
computer can do in the concept of a universal Turing machine, to
establish limits to what a computer can do in his solution to the {\it
  Entscheidungsproblem}, and yet to advocate nonetheless that
computers might imitate all intelligent activity. It remains a
challenging research program: to show that one limited branch of
mathematics, computation, might stand for all mathematical activity.

In the century since Turing's birth, the computer has become so
ubiquitous and the idea of computer as brain so commonplace that it
bears repeating that we must still think very long and hard about how
to construct a computer that can imitate a living, thinking
mathematician.

Proof assistant technology is still under development in labs; far
more is needed before it finds widespread adoption.  Ask any proof
assistant researcher, and you will get a sizable list of features to
implement: more automation, better libraries, and better user
interfaces!  Wiedijk discusses ten design questions for the next
generation of proof assistants, including the type system, which
axiomatic foundations of mathematics to use, and the language of proof
scripts~\cite{Wie10}.

Everyone actively involved in proof formalization experiences the
incessant barrage of problems that have been solved multiple times
before and that other users will have to solve multiple times again,
because the solutions are not systematic.  To counter this, the DRY
``Don't Repeat Yourself'' principle of programming, formulated
in~\cite{PP00}, has been carried to a refreshing extreme by Carette in
his proof assistant design.  For example, in his designs, a morphism
is defined only once, eliminating the need for separate definitions of
a morphism of modules, of algebras, of varieties, and so forth.
Carette's other design maxims include ``math has a lot of structure;
use it'' and ``abstract mathematical structures produce the best
code''~\cite{Car28p}.  Indeed, mathematicians turn to abstraction to
bring out relevant structure. This applies to computer code and
mathematical reasoning alike.  American Math Society guidelines for
mathematical writing apply directly to the computer: ``omit any
computation which is routine. \dots Merely indicate the starting
point, describe the procedure, and state the outcome''~\cite{DCFPS}
(except that computations should be automated rather than entirely
omitted).

We need to separate the concerns of construction, maintenance, and
presentation of proofs.  The construction of formal proofs from a
mathematical text is an extremely arduous process, and yet I often
hear proposals that would increase the labor needed to formalize a
proof, backed by secondary goals such as ease of maintenance, elegance
of presentation, fidelity to printed texts, and pedagogy.\footnote{%
  To explain the concept of {\it separation of concerns}, Dijkstra
  tells the story of an old initiative to create a new programming
  language that failed miserably because the designers felt that the
  new language had to look just like {\it FORTRAN} to gain broad
  acceptance. ``The proper technique is clearly to postpone the
  concerns for general acceptance until you have reached a result of
  such a quality that it deserves acceptance''
  \cite{DijkST}.} 
Better to avail ourselves of automation that was not available in the
day of paper proofs, and to create new mathematical styles suited to
the medium, with proofs that variously look like a
computer-aided design session, a functional program, or a list of
hypotheses as messages in gmail.
The most pressing concern is to 
reduce the skilled labor it takes a user to construct a formal proof
from a pristine mathematical text.

The other concerns of proof transformation should be spun off as separate
research activities: refactored proofs,
proof scripts optimized for execution time, translations into other
proof assistants, natural language translations, natural language
abstracts, probabilistically checkable proofs, searchable metadata
extracts, and proof mining.

\bigskip

For a long time, proof formalization technology was unable to advance
beyond the mathematics of the $19$th century, picking classical gems
such as the Jordan curve theorem, the prime number theorem, or
Dirichlet's theorem on primes in arithmetic progressions.  With the
Feit-Thompson theorem, formalization has risen to a new level, by
taking on the work of a Fields medalist.

At this level, there is an abundant supply of mathematical theorems to
choose from.  A Dutch research agenda 
lists the formalization of Fermat's Last Theorem as the first in a
list of ``Ten Challenging Research Problems for Computer
Science''~\cite{Berg}.  Hesselink predicts that this one
formalization project alone will take about ``fifty years, with a very
wide margin.''
Small pieces of the proof of Fermat, such as class field theory, the
Langlands-Tunnell theorem, or the arithmetic theory of elliptic curves
would be a fitting starting point.  The aim is to develop technologies
until formal verification of theorems becomes routine at the level of
Atiyah-Singer index theorem, Perelman's proof of the Poincar\'e
conjecture, the Green-Tao theorem on primes in arithmetic progression,
or Ng\^o's proof of the fundamental lemma.

% Hesselink quote: http://www.cs.rug.nl/~wim/fermat/wilesEnglish.html 

\bigskip

Starting from the early days of Newell, Shaw, and Simon's experiments,
researchers have dreamed of a general-purpose mechanical problem
solver.  Generations later, after untold trials, it remains an
unwavering dream.  I will end this section with one of the many
proposals for a general problem solving algorithm.  Kurzweil breaks
general problem solving into three phases:
\begin{enumerate} 
\item State your problem in precise terms.
\item Map out the contours of the solution space by traversing it
  recursively, within the limits of available computational resources.
\item Unleash an evolutionary algorithm to configure a neural net to
  tackle the remaining leaves of the tree.
\end{enumerate}
He concludes, ``And if all of this doesn't work, then you have a
difficult problem indeed''~\cite{Ku99}.  Yes, indeed we do!  Some day,
energy and persistence will conquer.

\newpage
\section{Issues of Trust}\label{sec:trust}

We all have first-hand experience of the bugs and glitches of
software.  We exchange stories when computers run amok.  Science
recently reported the story of a textbook ``The Making of a Fly'' that
was on sale at Amazon for more than 23 million dollars~\cite{Sci11}.
The skyrocketing price was triggered by an automated bidding war
between two sellers, who let their algorithms run unsupervised.  The
textbook's author, Berkeley professor Peter Lawrence, said he hoped
that the price would reach ``a billion.''
An overpriced textbook on the fly is harmless, except for students who
have it as a required text.  

But what about the Flash Crash on Wall
Street that brought a 600 point plunge in the Dow Jones in just 5
minutes at 2:41 pm on May 6, 2010?  According to the New York Times
\cite{NYT2010}, the flash crash started when a mutual fund used a
computer algorithm ``to sell \$4.1 billion in futures contracts.''
The algorithm was designed to sell ``without regard to price or
time.\dots [A]s the computers of the high-frequency traders traded
[futures] contracts back and forth, a `hot potato' effect was
created.''  When computerized traders backed away from the unstable
markets, share prices of major companies fluctuated even more
wildly. ``Over 20,000 trades across more than 300 securities were
executed at prices more than 60\% away from their values just moments
before'' \cite{SEC2010}. Throughout the crash, computers followed
algorithms to a T, to the havoc of the global economy.

% http://en.wikipedia.org/wiki/2010_Flash_Crash
% 2:41-2:46pm time from the NYT graphic.

\subsection{mathematical error}

Why use computers to verify mathematics?  The simple answer is that
carefully implemented proof checkers make fewer errors than
mathematicians (except J.-P. Serre).

Incorrect proofs of correct statements are so abundant that they are
impossible to catalogue.  Ralph Boas, former executive editor of Math
Reviews, once remarked that proofs are wrong ``half the
time''~\cite{Aus}. Kempe's claimed proof of the four-color theorem
stood for more than a decade before Heawood refuted
it~\cite[p.~115]{Mac}.  ``More than a thousand false proofs [of
Fermat's Last Theorem] were published between 1908 and 1912
alone''~\cite{Corry}.  Many published theorems are like the hanging
chad ballots of the 2000 U.S. presidential election, with scrawls too
ambivalent for a clear yea or nay.  One mathematician even proposed to
me that a new journal is needed that unlike the others only publishes
reliable results.  Euclid gave us a method, but even he erred in the
proof of the very first proposition of the Elements when he assumed
without proof that two circles, each passing through the other's
center, must intersect.  The concept that is needed to repair the gap
in Euclid's reasoning is an intermediate value theorem.  This defect
was not remedied until Hilbert's `Foundations of Geometry.'

Examples of widely accepted proofs of false or unprovable statements
show that our methods of proof-checking are far from perfect.
Lagrange thought he had a proof of the parallel postulate, but had
enough doubt in his argument to withhold it from publication.  In some
cases, entire schools have become sloppy, such as the Italian school
of algebraic geometry or real analysis before the revolution in rigor
towards the end of the nineteenth century.  Plemelj's 1908 accepted
solution to Hilbert's 21st problem on the monodromy of linear
differential equations was refuted in 1989 by Bolibruch.  Auslander
gives the example of a theorem\footnote{The claim was that every
  homogeneous plane continuum is a simple closed curve.}  published by
Waraskiewicz in 1937, generalized by Choquet in 1944, then refuted
with a counterexample by Bing in 1948~\cite{Aus}.  Another
example
%\footnote{I would like to thank J. Manfredi for this example.}
is the approximation problem for Sobolev maps between two
manifolds~\cite{Bethuel}, which contains a faulty proof of an
incorrect statement.  The corrected theorem appears in \cite{Hang}.
Such examples are so plentiful that a Wiki page has been set up to
classify them, with references to longer discussions at Math
Overflow~\cite{WikiPIP},~\cite{Over2},~\cite{Over1}.

% Plemelj: According to math overflow, ``Details are in the book The
% Riemann-Hilbert Problem by Anosov and Bolibruch (Vieweg-Teubner
% 1994), and a nice popular recounting of the story is in Ben
% Yandell's The Honors Class (A K Peters 2002).''

\begin{figure}[h!]
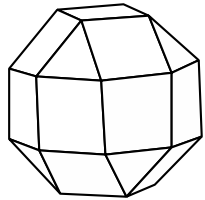

  \centering
\figPSEUDO
\caption{Throughout history, the pseudo rhombic cuboctahedron has been
  overlooked or misclassified.}
\label{fig:pseudo}
\end{figure}

Theorems that are calculations or enumerations are especially prone to
error.  Feynman laments, ``I don't notice in the morass of things that something, a
little limit or sign, goes wrong.\dots I have mathematically proven to myself
so many things that aren't true''
\cite[p.~885]{FeCo}. Elsewhere, Feynman describes two teams of
physicists who carried out a two-year calculation of the electron
magnetic moment and independently arrived at the same predicted value.
When experiment disagreed with prediction, the discrepancy was
eventually traced to an arithmetic error made by the physicists, whose
calculations were not so independent as originally
believed~\cite[p.~117]{FQED}.  Pontryagin and Rokhlin erred in
computing stable homotopy groups of spheres.  Little's tables of knots
from 1885 contains duplicate entries that went undetected until 1974.
In enumerative geometry, in 1848, Steiner counted $7776$ plane conics
tangent to $5$ general plane conics, when there are actually only
$3264$.  One of the most persistent blunders in the history of
mathematics has been the misclassification (or misdefinition) of
convex Archimedean polyhedra.  Time and again, the pseudo rhombic
cuboctahedron has been overlooked or illogically excluded from the
classification (Figure~\ref{fig:pseudo})~\cite{Gr11}.

% graphic by Hales

% references:
%  Pontryagin in overflow, Rokhlin in wiki.
%  Pontryagin in Ioan Mackenzie James, History of Topology p 567.
%  http://mathoverflow.net/questions/879 , calc examples.
%  http://mathoverflow.net/questions/35468
%  http://mathoverflow.net/questions/879
%  Lagrange in http://mathdl.maa.org/images/upload_library/22/Ford/Grabiner3-18.pdf

\subsection{In HOL Light we trust}

To what extent can we trust theorems certified by a proof assistant
such as HOL Light?  There are  various aspects to this question.
Is the underlying logic of the system consistent?  Are there any
programming errors in the implementation of the system?  Can a devious
user find ways to create bogus theorems that circumvent logic?  Are
the underlying compilers, operating system, and hardware reliable?

As mentioned above, formal methods represent the best cumulative effort of
logicians, computer scientists and mathematicians over the decades and
even over the centuries to create a trustworthy foundation for the
practice of mathematics, and by extension, the practice of science and
engineering.  

\subsection{a network of mutual verification}

John Harrison repeats the classical question {\it ``Quis custodiet
  ipsos custodes''} -- who guards the guards~\cite{HaSelf}?  How do we
prove the correctness of the prover itself?  In that article, he proves
the consistency of the HOL Light logic and the correctness of its
implementation in computer code.  He makes this verification in HOL
Light itself!  To skirt G\"odel's theorem, which implies that HOL
Light -- if consistent -- cannot prove its own consistency, he gives
two versions of his proof.  The first uses HOL Light to verify a
weakened version of HOL Light that does not have the axiom of
infinity.  The second uses a HOL Light with a strengthened axiom of
infinity to verify standard HOL Light.

Recently, Adams has implemented a version of HOL called HOL Zero.  His
system has the ability to import mechanically proofs that were
developed in HOL Light~\cite{Adams}.  He imported the
self-verification of HOL Light, to obtain an external verification.
You see where this is going.  As mechanical translation capabilities
are developed for proof assistants, it becomes possible for different
proof assistants to share consistency proofs, similar to the way that
different axiomatic systems give relative consistency proofs of one
another.  We are headed in the direction of knowing that if the logic
or implementation of one proof assistant has an error, then all other
major proof assistants must fail in tandem. Other self-verification
projects are Coq in Coq (Coc) and ACL2 in ACL2 (Milawa) ~\cite{Bar98},
~\cite{Dav09}.

\subsection{hacking HOL}

Of course, every formal verification project is a verification of an
abstract model of the computer code, the computer language, and its
semantics.  In practice, there are gaps between the abstract
model and implementation.

This leaves open the possibility that a hacker might find ways to
create an unauthorized theorem; that is, a theorem generated by some
means other than the rules of inference of HOL Logic.  Indeed, there
are small openings that a hacker can exploit.\footnote{For example,
  strings are mutable in HOL Light's source language, Objective CAML,
  allowing theorems to be maliciously altered.  Also, Objective CAML
  has {\it object magic}, which is a way to defeat the type system.
  These vulnerabilities and all other vulnerabilities that I know
  would be detected during translation of the proof from HOL Light to
  HOL Zero.  A stricter standard is Pollack consistency, which
  requires a proof assistant to avoid the appearance of
  inconsistency~\cite{Adams},~\cite{wiedijk:pollack}.  For example,
  some proof assistants allow the substitution of a variable
  whose name is a meaningless sequence of characters {\tt
    `n<0~$\land$~0' } for $t$ in $\exists n.~t<n$ to obtain a
  Pollack-inconsistency $\exists n.~{ {n<0~\land~0}}<n$.  } Adams maintains
a webpage of known vulnerabilities in his system and offers a cash
bounty to anyone who uncovers a new vulnerability.

%% HOL Light:
% let pollack_inconsistency = 
%  let thm1 = prove(  `?(n:num). (t < n)`,EXISTS_TAC`t + 1` THEN ARITH_TAC) in
%  let t = mk_var("n < 0 /\\ 0",`:num`) in
%    INST [(t,`t:num`)] thm1;;
%%

These documented vulnerabilities need to be kept in perspective.  They
lie at the fringe of the most reliable software products ever
designed. Proof assistants are used to verify the correctness of chips
and microcode~\cite{FoxArm6}, operating system kernels~\cite{seL4}, 
compilers~\cite{CC}, safety-critical software such as aircraft
guidance systems, security protocols, and mathematical theorems that
defeat the usual refereeing process.  

Some take the view that nothing short of absolute certainty in
mathematics gives an adequate basis for science.  Poincar\'e was less
exacting\footnote{``Il est donc inutile de demander au calcul plus de
  pr\'ecision qu'aux observations; mais on ne doit pas non plus lui en
  demander moins''~\cite{HPMC}.},
% quote appears in the introduction of the book.
only demanding the imprecision of calculation not to exceed
 experimental error.  As Harrison reminds us, 
``a foundational death spiral adds little value''~\cite{harrison-pm}.

\subsection{soft errors}\label{sec:soft}

Mathematicians often bring up the ``cosmic ray argument'' against the
use of computers in math.  Let's look at the underlying science.

A soft error in a computer is a transient error that cannot be
attributed to permanent hardware defects nor to bugs in software.
Hard errors -- errors that can be attributed to a lasting hardware
failure -- also occur, but at rates that are ten times smaller than
soft errors~\cite{MW04}.
Soft errors come from many sources. A typical soft error is caused by
cosmic rays, or rather by the shower of energetic neutrons they
produce through interactions in the earth's atmosphere.  A nucleus of
an atom in the hardware can capture one of these energetic neutrons
and throw off an alpha particle, which strikes a memory circuit and
changes the value stored in memory.  To the end user, a soft error
appears as a gremlin, a seemingly inexplicable random error that
disappears when the computer is rebooted and the program runs again.

As an example, we will calculate the expected number of soft errors in
one of the mathematical calculations of Section~\ref{sec:e8}.  The
Atlas Project calculation of the $E_8$ character table was a $77$ hour
calculation that required $64$ gigabytes RAM~\cite{AtlasSlides}.  Soft
errors rates are generally measured in units of failures-in-time
(FIT). One FIT is defined as one error per $10^9$ hours of operation.
If we assume a soft error rate of $10^3$ FIT per Mbit, (which is a
typical rate for a modern memory device operating at sea
level\footnote{The soft error rate is remarkably sensitive to
  elevation; a calculation in Denver produces about three times more
  soft errors than the same calculation on identical hardware in Boston.}
\cite{WP}),
 then we would expect there to be about $40$ soft
errors in memory during the calculation:
\[
\frac{10^3 \text{~FIT}}{1\text{~Mbit}} \cdot 64 \text{~GB} \cdot 77\text{~hours} =
\frac{10^3 \text{~errors~}}{10^9\text{~hours~}\text{Mbit}} \cdot
({64\cdot 8\cdot 10^3 \text{~Mbit}}) \cdot 77\text{~hours~} 
\approx 39.4 \text{~errors}.
\]
This example shows that soft errors can be a realistic concern in
mathematical calculations.  (As added confirmation, the $E_8$
calculation has now been repeated about $5$ times with identical
results.)

% references:
%The soft error rate of a memory device is typically in the range 1000
%to 5000 FIT per Mbit \cite{WP}.
% http://www.tezzaron.com/about/papers/soft_errors_1_1_secure.pdf
% Denver/Boston, 3 x in Mastipuram and Wee, EDN article.

In software that has been thoroughly debugged, soft errors become the
most significant source of error in computation.  Although there are
numerous ways to protect against soft errors with methods such as repeated calculations and error-correcting
codes, hardware redesign carries an economic cost.  In fact, soft errors are on
the rise through miniaturization: a smaller circuit generally has a lower
capacitance and responds to less energetic alpha particles than a larger
circuit.

Soft errors are depressing news in the ultra-reliable world of proof
assistants.  Alpha particles rain on perfect and imperfect software
alike.  In fact, because the number of soft errors is proportional to
the execution time of a calculation, by being slow and methodical, the
probability of a soft error during a calculation inside a proof
assistant can be much higher than the probability when done outside.

Soft errors and susceptibility to hacking have come to be more than a
nuisance to me.  They alter my philosophical views of the foundations
of mathematics.  I am a computational formalist -- a formalist who
admits physical limits to the reliability of any verification process,
whether by hand or machine.  These limits taint even the simplest
theorems, such as our ability to verify that $1+1=2$ is a consequence
of a set of axioms.  One rogue alpha particle brings all my schemes of
perfection to nought.  The rigor of mathematics and the reliability of
technology are mutually dependent; math to provide ever more accurate
models of science, and technology to provide ever more reliable
execution of mathematical proofs.

\newpage
\section{Concluding Remarks}

To everyone who has made it this far in this essay, I highly recommend
MacKenzie's book \cite{Mac}.  It written by a sociologist with a fine
sensitivity to mathematics.  The author received the Robert K. Merton
Award of the American Sociological Association in 2003 for this book.

  A few years ago, a special issue of the Notices of the
AMS presented a general introduction to formal
proofs~\cite{Hales:2008:formal},~\cite{Harrison:2008:formal},
~\cite{gonthier:2008:formal}, ~\cite{Wiedijk:2008:formal}.  I also
particularly recommend the body of research articles by Harrison,
Gonthier, and Carette.

\bigskip

I thank Adams (both Jeff and Mark), Urban, Carette, Kapulkin, Harrison, and Manfredi for conversations
about ideas in this article.

\raggedright
\bibliographystyle{amsalpha} % was plain %plainnat
\bibliography{/Users/thomashales/Desktop/googlecode/flyspeck/latex/bibliography/all}

\bigskip
\noindent
%\svninfo
\smallskip

\noindent

\end{document}